
\documentstyle{amsppt}

\TagsOnLeft
\topmatter
\title Complexity of Weakly Null Sequences
\endtitle
\author Dale E. Alspach
and Spiros Argyros \endauthor
\address Department of Mathematics,  Oklahoma State University,
Stillwater, OK 74078-0613
\endaddress
\address Department of Mathematics, University of Crete,
Heraklion, Crete\endaddress
\thanks{Research of the first author was
supported in part by NSF grant DMS-8602395.}\endthanks
\keywords{ordinal, index, complexity, unconditional, $L_1$-predual,
$C(K)$space}\endkeywords
\subjclass{46B20}\endsubjclass

\abstract{We introduce an ordinal index which measures the complexity
of a weakly null sequence, and show that a construction due to J.
Schreier can be iterated to produce for each $\alpha <\omega_1$, a
weakly null
sequence $\left(x^{\alpha}_n\right)_n$ in $C\left(\omega^{\omega^{
\alpha}}\right)$ with complexity $\alpha$.  As in the Schreier
example each of these is a sequence of indicator functions which is a
suppression-1 unconditional basic sequence.  These sequences are used
to construct Tsirelson-like spaces of large index.  We also show that
this new ordinal index is related to the Lavrentiev index of a Baire-1
function and use the index to sharpen some results of Alspach and
Odell on averaging weakly null sequences.}\endabstract

\endtopmatter
\document

\head 0.  Introduction\endhead

In this paper we investigate the oscillatory behavior of pointwise
converging sequences.  Our main tool is a new ordinal index which
measures the oscillation of such a sequence.  We show that there are
weakly null sequences of indicator functions in $C(K)$ with arbitrarily
large oscillation index and that the oscillation index is smaller than
other similar ordinal indices.  In particular the oscillation index of
a
pointwise converging sequence is directly compared to the Lavrentiev
index of its limit, the $\ell^1$-index defined by Bourgain, and the
averaging
index.  Many of the results are directly related to those in [A-O]
where the averaging index is used and we extend some results of
that paper.

The first example of a weakly null sequence with no subsequence
with averages going to zero in norm was constructed by Schreier
[Sch].  His example is a sequence of indicator functions and, as
observed by Pelczynski and Szlenk [P-S], can be realized on the space
of ordinals less than or equal to $\omega^{\omega}$ in the order
topology.  Because
these are indicator functions the failure of the Banach-Saks property
is solely dependent on the intersection properties of the sets.  Our
construction is also based on intersection properties of sets.  Thus
one viewpoint on the constructions in this paper is that there are
families of sets with very complicated intersection properties.  The
purpose of the ordinal index is to measure the complexity of these
intersection properties.  The examples of weakly null sequences that
we construct, like Schreier's example, are suppression-1 unconditional
basic sequences and can be considered as generalizations of Schreier's
construction.

Let us note that there is a strong relationship between this work
and some unpublished results of Rosenthal on the unconditional basic
sequence problem.  Rosenthal showed that any weakly null sequence
of indicator functions in $C(K)$ has a subsequence which is an
unconditional basic sequence.  In our work we found that Rosenthal's
notion of weakly independent sets fits naturally into our viewpoint
and we have incorporated some of Rosenthal's work into the
exposition.  While we believe that our work explains some difficulties
with unconditionality, it is based on intersection properties which
cannot provide a complete explanation.  For example, in contrast to
Rosenthal's result on sequences of indicator functions, it is known
that given a weakly null sequence one cannot always find a
subsequence which is an unconditional basic sequence [M-R] and [O].
Thus understanding the unconditional basic sequence problem is more
complicated than just understanding the oscillation properties of
weakly null sequences.  On the other hand it is clear that the
oscillation properties do play a fundamental role in unconditionality.
It would be interesting if there were some way to incorporate these
more subtle properties of weakly null sequences into an ordinal
index.

The paper is organized into six sections.  In the first we recall the
definitions of some ordinal indices and trees.  In the second we
introduce the oscillation index and an essentially equivalent index
which we call the spreading model index.  In the third we prove that
the oscillation index is essentially smaller than the $\ell^1$-index
and show
that it is related to the Lavrentiev index of a Baire-1 function.  In
the fourth section we define for each countable ordinal $\alpha$ a
weakly
null sequence of indicator functions and compute the oscillation index
of the sequence and the size of the smallest $C(K)$ space which can
contain the sequence.  In the fifth section we use an idea of Odell to
show that the construction in Section 4 can be used to construct
reflexive spaces similar to Tsirelson space with large oscillation
index.  In the sixth section we show that the averaging index, which
has a definition similar to that of the spreading model index, is not
smaller than the $\ell^1$ index.  In particular a space is constructed
which
does not contain $\ell^1$ but has averaging index $\omega_1$.  We also
extend some
results from [A-O] in order to better characterize those sequences
which can be averaged a predictable finite number of times in order
to get a weakly null sequence.

We will use standard terminology and notation from Banach space
theory as may be found in the books of Lindenstrauss and Tzafriri
[L-T,I] and [L-T,II] and Diestel [D1] and [D2].  If $\alpha$ is an
ordinal, we
will use $\alpha$ rather than $\alpha +1$ to denote the space of
ordinals less than
or equal to $\alpha$ in the order topology and $C(\alpha )$ to denote
the space of
continuous functions on $\alpha$.

\head 1.  Preliminaries\endhead

In this section we will recall the definitions of the $\ell^1$-index of
Bourgain, the Szlenk index, and the averaging index.  In order to
define the Bourgain $\ell^1$-index we first need to define trees and
some
related notions.
\definition{Definition}  Given a set $S$ a tree $\Cal T$ on $S$ is
subset of
$\cup^{\infty}_{n=1}S^n\cup S^{\infty}$ such that if $b\in \Cal T$ then
any initial
segment of $b$ is also in $\Cal T$, i.e., if $b\in \Cal T$ and $b=\left
(s_1,s_2,\ldots ,s_n,s_{n+1}\right)$ or
$b=\left(s_1,s_2,\ldots ,s_n,s_{n+1},\ldots\right)$, then $\left(s_
1,s_2,\ldots ,s_n\right)\in \Cal T$.  We will call the
elements of the tree {\it nodes\/} and say that the node $b$, as above,
is a
{\it successor\/} of or is {\it below $\left(s_1,s_2,\ldots ,s_n\right
)$\/}.  If $b$ and $c$ are nodes, we will
write $b>c$ to indicate that $b$ is below $c$.  In particular
$\left(s_1,s_2,\ldots ,s_n,s_{n+1}\right)$ is an {\it immediate
successor\/} of $\left
(s_1,s_2,\ldots ,s_n\right)$.  If a
node $x\in S^n$ then $n$ is the {\it length\/} or {\it level\/} of $
x$.  A {\it branch\/} of a tree is
a maximal linearly ordered subset of the tree under this natural
initial segment ordering.  A {\it subtree\/} of a tree $\Cal T$ is a
subset of $
\Cal T$
which is a tree.  A tree $\Cal T$ is {\it finitely branching\/} if the
number of
immediate successor nodes of any node is finite.  It is {\it dyadic\/}
if this
number is at most two for all nodes.  Finally a tree is {\it
well-founded \/}
if all of its branches are of finite cardinality.\enddefinition

For well-founded trees there is a standard way to define the ordinal
index of the tree.
\definition{Definition}  Suppose that $\Cal T$ is a well-founded tree
on a set $
X$.  Let
$\Cal T_0=\Cal T$ and for each ordinal $\alpha$ define
$$\Cal T^{\alpha +1}=\cup^{\infty}_{n=1}\big\{\left(x_1,x_2,\ldots
,x_n\right)\in \Cal T^{\alpha}:\text{ there is an }x\in X\text{ such
that}$$
$$\left(x_1,x_2,\ldots ,x_n,x\right)\in \Cal T^{\alpha}\big\}.$$
For a limit ordinal $\alpha$ let $\displaystyle{\Cal T^{\alpha}=\cap_{
\beta <\alpha}\Cal T^{\beta}}$.  Let $o(\Cal T)$ be the
smallest ordinal $\gamma$ such that $\Cal T^{\gamma}=\emptyset$.  This
is the order of the tree $
\Cal T$.
\enddefinition

Note that if $\overline x=\left(x_1,x_2,\ldots ,x_j\right)\in\!\Cal T^{
\alpha}\backslash \Cal T^{\alpha +1}$, the tree
$$\Cal T_{\overline x}=\cup^{\infty}_{n=1}\left\{\left(y_1,y_2,\ldots
,y_n\right):\overline x+\left(y_1,y_2,\ldots ,y_n\right)\right\}.$$
where $``+$'' indicates concatenation, is of order $\alpha$.

Now we are ready to define the $\ell^1$-index of Bourgain.
\definition{Definition}  Let $X$ be a Banach space and $\delta >0$.
Let
$$\align\!\Cal T(X,\delta )=\cup^{\infty}_{n=1}\big\{\big(x_1,x_2,&
\ldots ,x_n\big)\in X^n:||x_j||\le 1\text{ for all }j\le n\text{ and
}\\
&||\sum_j\lambda_jx_j||\ge\delta\sum |\lambda_j|\text{ for all
}(\lambda_
j)\in \Bbb R^n\big\}.\endalign$$
The $\delta\text{ }\ell^1$-index of $X$ is $\ell(X,\delta )=o(\Cal T
(X,\delta ))$.  If $\Cal T(X,\delta )^{\alpha}\neq\phi$ for all
countable $\alpha$ then $\ell (X,\delta )=\omega_1$.\enddefinition
\vskip .15in
It is easy to see that if $\ell^1$ is isomorphic to a subspace of $
X$ then
$\ell (X,\delta )=\omega_1$ for some $\delta >0$ and thus by a result
of James [J] for all
$\delta$, $0<\delta <1$.  The converse is also true but requires the
fact that if
$\Cal T^{\alpha}(X,\delta )\neq\emptyset$ for all countable $\alpha$
then the tree has an infinite branch.
(See [Bo].)  Bourgain observed that the examples employed by Szlenk
[Sz] to show that there are separable reflexive spaces with Szlenk
index greater than any given countable ordinal also have large
$\ell^1$-index.  Note that this indicates that the Baire-1 functions in
$
X^{**}$
provide only a very weak indication as to the $\ell^1$-index of $X$.

Next we will define Bourgain's index for Boolean independence.  This
index is a technical convenience which we will use in establishing
lower bounds on the $\ell^1$-index.
\definition{Definition}  Let $K$ be a set and $\left\{\left(A_n,B_
n\right)\right\}$ be a sequence of pairs of
subsets of $K$.  Let
\TagsOnLeft

$$\align \Cal T(\{(A_n,B_n)\})=\cup^{\infty}_{k=1}\big\{(&n_1,n_2,
\ldots ,n_k)\in \Bbb N^k:\\
&\text{for all }(\epsilon_1,\epsilon_2,\ldots ,\epsilon_k)\in \{-1
,1\}^k,\cap^k_{i=1}\epsilon_iA_{n_i}\neq\emptyset\big\},\endalign$$
where $(-1)A_i=B_i$.  Let $\Cal B(\{(A_n,B_n)\})$ be the order of the
tree
$\Cal T(\{(A_n,B_n)\})$.
\enddefinition

We will usually take pairs of disjoint sets $(A_n,B_n)$ when we employ
this index.  Note that by Rosenthal's theorem on $\ell^1$, [R], a
bounded
sequence $(x_n)$ in a Banach space $X$ has a subsequence equivalent to
the unit vector basis of $\ell^1$ if (and only if)  it has no weak
Cauchy
subsequence.  Moreover for any such sequence there is a subsequence
$(x_n)_{n\in M}$ and two numbers $\delta >0$ and $r$ such that if
$$A_n=\left\{x^{*}\in X^{*}:x^{*}(x_n)\ge r+\delta ,||x^{*}||\le
1\right
\}$$
and
$$B_n=\left\{x^{*}\in X^{*}:x^{*}(x_n)\le r,||x^{*}||\le 1\right\}$$
then $\Cal T\left(\left\{\left(A_n,B_n\right):n\in M\right\}\right
)=\cup^{\infty}_{k=1}M^k$ and thus
$\Cal B$ $\left(\left\{\left(A_n,B_n\right):n\in \Bbb N\right\}\right
)=\omega_1$.  Hence if we consider such sets $(A_n,B_n)$,
the index is $\omega_1$ for some $\delta >0$ and $r$ if and only if $
(x_n)$ has a
subsequence with no weak Cauchy subsequence.  In particular if $(x_
n)$
converges $w^{*}$ sequentially to some $x^{**}\in X^{**}$ the index is
countable.

Now we wish to recall some other ordinal indices which are more
classical in spirit.  Each of these indices is defined in terms of the
oscillation of sequences of functions on the unit ball of the dual in
the $w^{*}$ topology.  First we recall the Szlenk index.  Let $X$ be a
Banach space and for each $\alpha <\omega_1$ let
$$\align P_{\alpha +1}&\left(\varepsilon
,B_X,B_{X^{*}}\right)=\big\{x^{*}
:\text{ there exists }x_n^{*}\in P_{\alpha}\left(\epsilon ,B_X,B_{
X^{*}}\right)\text{ and}\\
&x_n\in B_X\text{ such that }x^{*}_n@>{w^{*}}>>x^{*},x_n@>{w}>>0,\text{
and }\lim
x^{*}_n(x_n)\ge \varepsilon\big\}\endalign$$
At a limit ordinal $\beta$ we take the intersection, that is,
$$P_{\beta}\left(\epsilon ,B_X,B_{X^{*}}\right)=\cap_{\alpha <\beta}
P_{\alpha}\left(\epsilon ,B_X,B_{X^{*}}\right).$$

Now fix $\{e_j\}$ a normalized weakly null sequence in $X$.  The next
two
ordinal indices are defined for each such sequence.  The first
occurred in a paper of Zalcwasser [Z] and in a paper of Gillespie and
Hurwitz [G-H] and was used to prove that a pointwise converging
sequence of bounded continuous functions on a compact metric space
has a sequence of convex combinations going to zero in norm.
Let
$$\align&Z_{\alpha +1}\left(\epsilon
,\left\{e_j\right\},B_{X^{*}}\right
)=\big\{x^{*}:\text{there exists }x^{*}_n\in Z_{\alpha}\left(\epsilon
,\left\{e_j\right\},B_{X^{*}}\right)\text{ and}\\
&e_{j_n}\text{ such that }x^{*}_n@>{w^{*}}>>x^{*}\text{ and }\lim
|x^{*}_n(e_{j_{2n}})-x^{*}_n(e_{j_{2n-1}})|\ge\epsilon\big\}.\endalign$$
As before $\displaystyle{Z_{\beta}=\cap_{\alpha <\beta}Z_{\alpha}}$,
for a limit ordinal $
\beta$.

In [A-O] the following index was introduced in order to obtain some
more precise information about the nature of the convex combinations
obtained in these early papers in the context of weakly null
sequences in a Banach space.

Let $A_{\alpha +1}\left(\epsilon ,\left\{e_j\right\},B_{X^{*}}\right
)=\big\{x^{*}:\text{for every neighborhood }\Cal N\text{ of }x^{*}$
$$\text{relative to }A_{\alpha}\left(\epsilon ,\left\{e_j\right\},
B_{X^{*}}\right)\text{there exists an infinite}$$
$$\text{set }L\subset \Bbb N\text{ such that }\ell^1-SP\big(\big\{
e_{i|\Cal N}\big\}_{i\in L}\big)\ge\epsilon\big\}$$
where
$$\align \ell^1-SP\left(\left\{e_i\right\}\right)=\\
&\lim_k\lim_m\inf\left\{k^{-1}\left|\left|\sum^k_{i=1}e_{n_i}\right
|\right|:m\le n_1<n_2<\ldots <n_k\right\}.\endalign$$
As before $\displaystyle{A_{\beta}=\cap_{\alpha <\beta}A_{\alpha}}$, if
$
\beta$ is a limit ordinal.  We
will refer to this as the averaging index.

In each case the successor set (if non-empty) is a $w^{*}$ closed
nowhere
dense subset of the set.  (For the case of the Szlenk sets we need to
assume that $X^{*}$ is separable.)  Thus the Baire Theorem gives in
each
case a largest ordinal $\alpha$ with the $\alpha$th set non-empty and
the $
(\alpha +1)$th
empty.  Denote $o(P,\epsilon )$, $o(Z,\epsilon )$, and $o(A,\epsilon
)$ as the indices for the Szlenk,
Zalcwasser, and average sets, respectively.  Clearly $o(P,\epsilon
)\ge o(Z,2\epsilon )$
and $o(Z,\epsilon )\ge o(A,\epsilon )$.

\head 2.  Weakly Null Sequences and the $\ell^1$-index\endhead

We wish to introduce an index similar to the Szlenk index which will
measure $\ell^1$-ness in a different way than the $\ell^1$-index of
Bourgain.
\definition{Definition}  Let $(K,d)$ be a Polish space and let $(f_
n)$ be a pointwise
convergent sequence of continuous functions on $K$.  Fix $\epsilon
>0$ and let
$$A_{n,m}^{+}=\left\{k\in K:f_n(k)-f_m(k)>\epsilon\right\}$$
and
$$A_{n,m}^{-}=\left\{k\in K:f_n(k)-f_m(k)<-\epsilon\right\},$$
For each countable ordinal $\alpha$ we define inductively a subset of $
K$ by
$$\Cal O^0\left(\epsilon ,\left(f_n\right),K\right)=K,$$
$$\align
\Cal O^{\alpha +1}\left(\epsilon ,\left(f_n\right),K\right)=&\big
\{k\in \Cal O^{\alpha}\left(\epsilon ,\left(f_n\right),K\right):\text{
for every neighborhood }\Cal N\text{ of }k\text{ there is an }\\
&N\in
\Bbb N
\text{ such that for all } n\ge N\text{ there
exists an }M\in \Bbb N\text{ such that }\\&
\cap_{m\ge M}A_{n,m}^{+}\cap \Cal O^{\alpha}\left(\epsilon ,\left
(f_n\right),K\right)\cap \Cal N\neq\emptyset\text{
or }
\\&\cap_{m\ge M}A_{n,m}^{-}\cap \Cal O^{\alpha}\left(\epsilon ,\left
(f_n\right),K\right)\cap \Cal N\neq\emptyset \}.\endalign$$
If $\beta$ is a limit ordinal,
$\displaystyle{\Cal O^{\beta}\left(\epsilon ,\left(f_n\right),K\right
)=}$$\displaystyle{\cap_{\alpha <\beta}\Cal O^{\alpha}\left(\epsilon
,\left(f_n\right
),K\right)}$.
\enddefinition

Note that if $(f_n)$ converges pointwise to 0 and $\epsilon^{\prime}
<\epsilon$, then for large
enough $M$
$$\left\{k:f_n(k)\ge\epsilon^{\prime}\right\}\supset\cap_{m\ge M}A_{
n,m}^{+}\supset\left\{k:f_n(k)\ge\epsilon\right\}$$
Thus $\Cal O^{\alpha}\left(\epsilon ,\left(f_n\right),K\right)$ is
essentially the Szlenk index set, $
P_{\alpha}\left(\epsilon ,\left(f_n\right),K\right)$,
except that we do not allow the use of subsequences.  Thus
$P_{\alpha}\left(\epsilon^{\prime},\left(f_n\right),K\right)\supset
\Cal O^{\alpha}\left(\epsilon ,\left(f_n\right),K\right)$ in this case
for all $
\epsilon^{\prime}<\epsilon$.  Also it is
easy to see that $\Cal O^{\alpha}\left(\epsilon ,\left(f_n\right),
K\right)$ is always closed.  If the limit $f$ is
not continuous, the relationship with the Szlenk sets is less clear.
However in the definition of the Szlenk sets it is possible to use
weak Cauchy sequences in place of weakly null sequences.  The sets
obtained in this way behave in essentially the same way as the
original Szlenk sets provided the dual is separable.  Thus if $[f_
n]^{*}$ is separable $\Cal O^{\alpha}\left(\epsilon ,\left(f_n\right
),K\right)$ will be empty
for some $\alpha <\omega_1$.  Actually the index is always countable
but this
will be a consequence of our result relating this index to the
$\ell^1$-index.

This discussion (and some later results) justifies the following.

\definition{Definition}  Suppose $(f_n)$ is a pointwise converging
sequence of
continuous functions on a Polish space $K$.  $\Cal O(\epsilon )=\Cal
O\left
(\epsilon ,\left(f_n\right),K\right)$ is the
largest ordinal $\alpha$ such that $\Cal O^{\alpha}\left(\epsilon
,\left(f_n\right),K\right)$ is non-empty.  We will
refer to $\Cal O(\epsilon )$ as the $\epsilon$ {\it oscillation
index\/} of the sequence $\left
(f_n\right)$.
(Some authors, e.g., [K-L], use the term oscillation index for an
ordinal index defined in terms of the oscillation of a fixed function
near a point.  Here we only apply the term to sequences so no
confusion should result.)

It will be convenient to use a more restrictive index at times.
Let\linebreak
$\Cal O_{+}^0\left(\epsilon ,\left(f_n\right),K\right)=K$.  For each
countable ordinal $
\alpha$ define inductively a
subset of $K$ by
$$\align \Cal O_{+}^{\alpha +1}\left(\epsilon ,\left(f_n\right),K\right
)&=\big\{k\in \Cal O_{+}^{\alpha}\left(\epsilon ,\left(f_n\right),
K\right):\text{ for every neighborhood $\Cal N$ of $k$}\\
&\text{there is an }N\in \Bbb N\text{ such that for all }n\ge N\text{
there exists }\\
&\text{an }M\in \Bbb N\text{ such that }\cap_{m\ge M}A_{n,m}^{+}\cap
\Cal O_{+}^{\alpha}\left(\epsilon ,\left(f_n\right),K\right)\cap \Cal N
\neq\emptyset\big\}.\endalign$$
If $\beta$ is a limit ordinal,
$\text{$\Cal O_{+}^{\beta}$ $\left(\epsilon ,\left(f_n\right),K\right
)=$ $\cap_{\alpha <\beta}\Cal O_{+}^{\alpha}$ $\left(\epsilon ,\left
(f_n\right),K\right)$}$.  As above the
positive $\epsilon$ oscillation index $\Cal O_{+}(\epsilon )$ is the
largest ordinal $
\alpha$ such that
$\Cal O_{+}^{\alpha}\left(\epsilon
,\left(f_n\right),K\right)\neq\emptyset$.  In a similar way we define $
\Cal O_{-}^{\alpha}\left(\epsilon ,\left(f_n\right),K\right)$ and
$\Cal O_{-}\left(\epsilon ,\left(f_n\right),K\right)$.\enddefinition

\remark{Remark 2.1}  We do not need to have the sequence converging
pointwise for the definition of the index to make sense.  However, if
$(f_n)$ has no pointwise convergent subsequence then, by Rosenthal's $
\ell^1$
theorem [R], there is an infinite set $M$ and real numbers $\delta
>0$ and $r$
such that $\left(\left\{k:f_n(k)\ge r+\delta\right\},
\left
\{k:f_n(k)\le r\right\}\right)_{n\in M}$ is Boolean independent.
From this it follows that there is a Cantor set $C\subset K$ so that
relative to $C$ and taking $n$ and $m$ from $M$, $\displaystyle{\cap_{
m>n}A_{n,m}}$
is an $\epsilon_n$-net in $C$ and $\epsilon_n\to 0$ as $n\to\infty$.
Hence $
\Cal O^1\left(\delta /2,\left(f_n\right)_{n\in M},C\right)=C$.
Thus for this subsequence the oscillation index is $\omega_1$.
\definition{Definition}  The $\epsilon${\it -oscillation index of a
Banach space\/} is the
supremum of the $\epsilon$-oscillation index over all weak Cauchy
sequences
in the unit ball of the space as functions on the dual ball with the
$w^{*}$ topology.\enddefinition
\vskip .15in
The next lemma follows easily from the definitions but we state it
for future reference.\endremark

\proclaim{Lemma 2.2}  If $K$ is a $w^{*}$ closed subset of the dual
ball of
$X$, then $\Cal O^{\alpha}\left(\epsilon ,\left(f_n\right),K\right
)\subset \Cal O^{\alpha}\left(\epsilon
,\left(f_n\right),B_{X^{*}}\right
)$ for all $\alpha$ and $\epsilon$.  In
particular if $T$ is a bounded operator from a Banach space $X$ to a
Banach space $Y$, then
$$T^{*}\Cal O^{\alpha}\left(\epsilon ,\left(Tf_n\right),B_{Y^{*}}\right
)=\Cal O^{\alpha}\left(\epsilon ,\left(f_n\right),T^{*}B_{Y^{*}}\right
)\subset \Cal O^{\alpha}\left(\epsilon /||T||,\left(f_n\right),B_{
X^{*}}\right),$$
and consequently
$$\Cal O\left(\epsilon ,\left(Tf_n\right),B_{Y^{*}}\right)\le \Cal
O\left
(\epsilon /||T||,\left(f_n\right),B_{X^{*}}\right).$$
\endproclaim

On the other hand we do not know if there is an essential
equivalence between\linebreak $\sup_M\Cal
O\left(\varepsilon,\left(f_n\right)_{
n\in M},B_{X^{*}}\right)$ and
$\sup_M\Cal O\left(\varepsilon\lambda,\left(f_n\right)_{n\in M},K\right
)$ if $K$ is a $\lambda$ norming subset of $B_{X^{*}}$.  Simple
examples show that considering subsequences is necessary in order
for the indices to be approximately the same size.

Next we will present some examples in which we can compute some
bounds on the oscillation index and compare this index to the Bourgain
$\ell^1$-index.  The first nontrivial example we wish to present is the
Schreier sequence, [Sch].

\example{Example 1}  Let $\Cal F_1=\left\{F\subset \Bbb N:\min F\ge
\text{card}F\right\}\cup \{\emptyset \}$.  We will identify\linebreak
$\Cal F_1$
with $\left\{1_F:F\in \Cal F_1\right\}$ and note that the later is a
countable compact
(metric) space in the topology of pointwise convergence on $\Bbb N$.
Define $f_n(F)=1_F(n)$.  Because each $F$ is finite $(f_n)$ converges
pointwise to 0 on $\Cal F_1$.  It is easy to see that each $f_n$ is
continuous on
$\Cal F_1$.  A computation shows that
$$\Cal O^j(1,(f_n),\Cal F_1)=\cup_{m>j}\left\{F\in \Cal F_1:\min F
\ge m\text{ and card }F\le m-j\right\}\cup \{\emptyset \}.$$
Thus $\Cal O^{\omega}(1,(f_n),K)=\{\emptyset \}$ and $\Cal
O(1)=\omega$.  Because these are indicator
functions the same is true for all $\epsilon ,$ $1\ge\epsilon >0$.
Note that this is the
maximal index because $\Cal F_1$ is homeomorphic to $\omega^{\omega}$
in the order
topology and hence there are only $\omega$ non-empty topological
derived
sets.

It follows from the combinatorial properties of $\Cal F_1$ that $[
f_n]$ contains
$\ell_n^{1\prime}s$ uniformly and thus
$\ell\left(\left[f_n\right]\right
)\ge\omega$.  (See [Pel-Sz] or Section 4.)
Also it is not hard to see that the $\epsilon$-oscillation index of $
C(\omega^{\omega})$ is
essentially the same as the Szlenk index, $\omega\left[\frac
1{\epsilon}\right
]$.\endexample
\example{Example 2}  Tsirelson space $T$.  (See [C-S] or Section 5.)
The natural
unit vector basis of $T$ is a weakly null sequence.  $T$ is reflexive
and
contains many $\ell_n^{1\prime}s$.  Thus $\ell (T,\delta )\ge\omega$,
for $
\delta =\frac 14$, for example.  The
computation of $\Cal O^j\left(\epsilon
,\left(e_n\right),B_{T^{*}}\right
)$ does not seem to be easy.  However
note that the functionals $1_F$ for $F\in \Cal F_1$ in example 1 above
act on $
T$
by $\displaystyle{1_f\left(\sum a_ne_n\right)=\sum_{n\in F}a_n}$ and
all $
||1_F||_{T^{*}}\le 2$.  Moreover
if $1_{F_k}\to 1_F$ pointwise on $\Bbb N$, then $1_{F_k}(x)\to 1_F
(x)$ for all $x\in T$.  Therefore
the map $S:T\to C\left(\!\Cal F_1\right)=C(\omega^{\omega})$ defined by
$
Sx(F)=1_F(x)$ is well defined
and bounded by 2 and $Te_n=f_n$.  It follows from Example 1 and
Lemma 2.2 that $\Cal O^{\omega}\left(\frac 12,\left(e_n\right),B_{
X^{*}}\right)\neq\emptyset$.\endexample
\example{Example 3}  $\ell^1$.  There are no non-trivial weakly null
sequences in $
\ell^1$
and thus the oscillation index is 0.  On the other hand the $\ell^
1$-index is
$\omega_1$.
\vskip .15in
The anomalous behavior of the index for $\ell^1$ can be corrected if we
allow general sequences and modify the definition of $\Cal O^{\alpha
+1}$.  However
this seems pointless in view of Rosenthal's characterization of $\ell^
1$.\endexample

\remark{Remark 2.3}  Suppose that $(f_n)$ is a pointwise converging
sequence of uniformly bounded continuous functions on a compact
metric space $K$.  Then the mapping $S$ from $K$ into $c$ (the space of
convergent sequences under the sup norm) $S(k)=(f_n(k))$ is continuous
if the range is given the $\sigma (c,\ell^1(\Bbb N))$ topology.  If the
functions
$(f_n)$
separate points on $K$ then $S$ is also one-to-one.  From this
viewpoint
we are investigating $\sigma (c,\ell^1(\Bbb N))$ compact subsets of the
ball of $
c$.  The
oscillation set $\Cal O^{\alpha}\left(\epsilon
,\left(f_n\right),K\right
)$ is mapped by $S$ to the set
$\Cal O^{\alpha}\left(\epsilon ,\left(e_n\right),S(K)\right)$, where $
e_n$ denotes the functional evaluation at $n$,
and thus the index may be computed in $c$.  If the sequence $(f_n)$
converges to 0 then the sets are actually in $c_0$ and the topology
$\sigma\left(c_0,\ell^1(\Bbb N)\right)$ is the weak topology.\endremark

\head 3.  Comparison with the $\ell^1$-index\endhead

The $\ell^1$-index defined by Bourgain [Bo] measures the degree to
which
$\ell^1$ isomorphically embeds in the space.  Bourgain showed that this
index is related to an ordinal index of Baire-1 functions in the second
dual.  In the next section we will show that the presence of certain
weakly null sequences can also raise the $\ell^1$-index.  Thus in
contrast
to the Baire-1 case the pointwise limit itself provides no information
but rather the sequence carries the information.  Our method is to
use the oscillation index defined in the previous section.  In this
section we will show that this oscillation index is essentially bounded
above by the $\ell^1$-index, and prove that any sequence of continuous
functions converging pointwise to a Baire-1 function of index $\alpha$
must
have a large oscillation index.  We then get Bourgain's result as a
corollary.  These ideas are also related to some work of Haydon,
Odell, and Rosenthal [H-O-R] on what they term Baire-1/2 and
Baire-1/4 functions in the second dual.

Now we will prove that a large oscillation index implies a large
$\ell^1$-index.  For an ordinal $\alpha =\omega^{\gamma}k+\beta$ with $
\beta <\omega^{\gamma}$ and $k\in \Bbb N$, we define
$\alpha /2=\omega^{\gamma}[(k+1)/2]$, where [ ] denotes the greatest
integer function.

\proclaim{Theorem 3.1}  If $(f_n)$ is a pointwise converging sequence
on
a compact metric \linebreak space $K$ and $\Cal
O^{\alpha}\left(\epsilon
,\left(f_n\right),K\right)\neq\emptyset$ then if
$\epsilon^{\prime}<\epsilon ,\text{
}\ell\left(\left[f_n\right],\epsilon^{
\prime}/2\right)\ge\alpha /2$.  Moreover there is an $\ell^1$-index
tree on
$(f_n)$ with index $\alpha /2$.\endproclaim

The proof is easier if we assume that $\Cal
O_{+}^{\alpha}\left(\epsilon
,\left(f_n\right),K\right)\neq\emptyset$, and in this
case we get $\alpha$ rather than $\alpha /2$ for the $\ell^1$-index.
The proof naturally
divides into two parts:  first a reduction to the case $\Cal O_{+}^{
\alpha /2}\neq\emptyset$ and
then the proof of the result in this case.  The reduction is proved as
Lemma 3.4.  The main idea in the remainder of the argument is to
construct a tree of Boolean independent pairs of sets with large
order where the pairs of sets are subsets of the $A_{n,m}^{+}$ $^{
\prime}s$.
Consequently we will actually construct our $\ell^1$ tree on
$\left(f_n-f_m\right)_{n,m\in \Bbb N}$ with constant $\epsilon /2$.  A
final argument is needed to get
a $\ell^1$ tree on $\left(f_n\right)$.

Before we proceed to the proof we need a few lemmas which
describe some sufficient conditions for a tree to be a $\ell^1$-tree.
The
following is an unpublished lemma of Rosenthal.  (We actually only
use the weaker version in which $r$ does not depend on
$n$.)

\proclaim{Lemma 3.2}   Suppose $\left(f_n\right)_{n=1}^m$ is a finite
sequence of norm
one functions on $K,\text{ }\delta >0$, and for each $n$ there is a
number $
r_n$ such
that if $A_n=\left\{f_n\ge r_n+\delta\right\}$ and $B_n=\left\{f_n
\le r_n\right\},\text{ }\left\{\left(A_n,B_n\right)\right\}$ is Boolean
independent.  Then $\left(f_n\right)$ is $2/\delta$ equivalent to the
unit vector basis of
$\ell_m^1$.\endproclaim

\demo{Proof}  Suppose that $a_n\in \Bbb R$ for all $n$. Let $F=\left
\{n:a_n\ge 0\right\}$ and
$G=\left\{n:a_n<0\right\}$.  Let $\displaystyle{t\in\cap_{
n\in F}A_n\cap\cap_{n\in G}B_n}$ and
$\displaystyle{t^{\prime}\in\cap_{n\in G}A_n\cap\cap_{n\in F}B_n}$.
Then
$$|\sum a_nf_n(t)-\sum a_nf_n(t^{\prime})|\ge\sum a_n[f_n(t)-f_n(t^{
\prime})]\ge\sum |a_n|\delta .$$
Indeed, if $n\in F$, $a_n\ge 0$, $f_n(t)\ge r_n+\delta$ and $f_n(t^{
\prime})\le r_n$ and if $n\in G$,
$a_n<0$, and $f_n(t)\le r_n$ and $f_n(t^{\prime})\ge r_n+\delta$.
Hence either
$|\sum a_nf_n(t)|\ge\sum |a_n|\delta /2$ or $|\sum a_nf_n(t^{\prime}
)|\ge\sum |a_n|\delta /2$.\qed\enddemo

We would like to construct a Boolean independent tree of pairs of
subsets of $K$ where the pairs are related to the elements of the
sequence $\left(f_n\right)$.  If we could, for example, choose sets of
the form $
A_n$
and $B_n$ as in the lemma above, we would immediately get a $\ell^
1$ tree
on $\left(f_n\right)$.  Actually somewhat less is sufficient.

\proclaim{Lemma 3.3}  Suppose that $\left(x_n\right)$ is a sequence of
functions
on a set $K$ with values in [-1,1] and $\Cal T$ is a tree on $\left
(x_n\right)$ of order
$\alpha <\omega_1$.
Further assume that there is a $\delta >0$ such that for each branch $
\Cal B$ of
$\Cal T$ there is a mapping $\rho_{\Cal B}:\Cal B\to 2^K\times 2^K$
such that

\itemitem{a)}  if $\rho_{\Cal B}(x_1,x_2,$ $\ldots ,x_n)$ $=(A,B)$ then
$
A\subset\left\{x_n\ge r+\delta\right\}$ and
$B\subset\left\{x_n\le r\right\}$ for some $r\in \Bbb R$.  ($r$ may
depend on $
n$ and $\Cal B$ $.$)
\itemitem{b)}  $\rho_{\Cal B}(\Cal B)$ is a Boolean independent
sequence of sets.\newline

Then $\Cal T$ is an $\ell^1$ tree of order $\alpha$ with constant $
\delta /2$.\endproclaim

\demo{Proof}  According to the previous lemma if
$\Cal B=\left\{\left(x_1\right),\left(x_1,x_2\right),\ldots\right\}$
then $\left
(x_i\right)$ is $2/\delta$ equivalent to the unit vector
basis of $\ell^1$.  Thus each branch of $\Cal T$ also satisfies the
requirements
for an $\ell^1$ tree, and thus we have an $\ell^1$-tree of order $
\alpha$.\qed\enddemo

Our next lemma allows us to reduce to the case of positive
oscillation index.  Below $p(\alpha )=\inf\left\{\beta +\rho :\rho
+\beta =\alpha\right\}$.  In particular if
$\alpha =\omega^{\gamma}k+\beta$, where $\beta <\omega^{\gamma}$, then
$
p(\alpha )=\omega^{\gamma}k$.

\proclaim{Lemma 3.4}  If $\Cal O^{\alpha}\left(\delta,\left(f_n\right
),K\right)\neq\emptyset$ and $\gamma\le\alpha /2$, then for
$\epsilon =+$ or $-$,\linebreak $\Cal
O^{\gamma}_{\epsilon}\left(\delta,\left
(f_{n_j}\right),K\right)\neq\emptyset$, for some subsequence $\left
(f_{n_j}\right)$.\endproclaim

\demo{Proof}  The proof is by induction on $\alpha$.  We will actually
prove that if\linebreak $t\in \Cal O^{\alpha}\left(\epsilon
,\left(f_n\right
),K\right)$ then there is an infinite set $L\subset \Bbb N$
and ordinals $\gamma$ and $\lambda$ such that
$$\min\left(\gamma +\lambda ,\lambda +\gamma\right)\ge p(\alpha )$$
and
$$t\in \Cal O_{+}^{\gamma}\left(\delta ,\left(f_n\right)_{n\in L}
,K\right)\cap \Cal O_{-}^{\lambda}\left(\delta ,\left(f_n\right)_{
n\in L},K\right).$$

Suppose that $p(\alpha )=\omega^v\cdot k$.  Because the result depends
only on $
p(\alpha )$,
we need only consider ordinals $\alpha$ of the form $\omega^v\cdot
k$ in the induction.

If $\alpha =1$, let $t\in \Cal O^{\alpha}\left(\delta ,\left(f_n\right
),K\right)$ and for $\epsilon =+$ or $-$ and $i\in \Bbb N$ consider
the set
$$N_{\epsilon}^i=\left\{n:\text{ there exists an }M\in N\text{ such
that }
\cap_{m\ge M}A_{nm}^{\epsilon}\cap \Cal N_i\neq\emptyset\right\}$$
where $\Cal N_i$ is a decreasing sequence of neighborhoods with
intersection
$\{t\}$.   For at least one choice of $\epsilon ,\text{ }N_{\epsilon}^
i$ is infinite for all $i$.  For that
$\epsilon$ let $L$ be an infinite subset of $\Bbb N$ such that
$L\backslash
N_{\epsilon}^i$ is finite for each $i$.
 Clearly $t\in \Cal O_{\epsilon}^1\left(\delta ,\left(f_n\right)_{
n\in L},K\right)$.

Now suppose that the lemma is true for all $\beta <\alpha$ and let
$t\in \Cal O^{\alpha}\left(\delta ,\left(f_n\right),K\right)$.  If $
\alpha =\omega^v$, let $\alpha_i\uparrow\alpha$.  The inductive
assumption
implies that there are sequences $\gamma_i$ and $\lambda_i$ and $L_
i\subset \Bbb N$ such that
$$\min\left(\gamma_i+\lambda_i,\gamma_i+\lambda_i\right)\ge p(\alpha_
i)$$
and
$$t\in \Cal O_{+}^{\gamma_i}\left(\delta,\left(f_n\right)_{n\in L_i}
,K\right)\cap \Cal O_{-}^{\lambda_i}\left(\delta,\left(f_n\right)_{
n\in L_i},K\right).$$
We may assume that $L_i\subset L_{i-1}$ for all $i$ and that the
sequences $\left
(\gamma_i\right)$
and $\left(\lambda_i\right)$ are non-decreasing.  It now follows that
if $
L$ is an infinite
subset of $\Bbb N$ such that $L\backslash L_i$ is finite for all $
i$ then
$$t\in \Cal O_{+}^{\gamma}\left(\delta,\left(f_n\right)_{n\in
L},K\right
)\cap \Cal O_{-}^{\lambda}\left(\delta,\left(f_n\right)_{n\in
L},K\right
),$$
where $\gamma =\lim \gamma_i$ and $\lambda =\lim \lambda_i$.  We have
that $\min\left
(\gamma +\lambda ,\lambda +\gamma\right)\ge p(\alpha_i)$
for all $i$ and $\min\left(\gamma +\lambda ,\lambda +\gamma\right)
\ge\lim p(\alpha_i)$.  Thus if $p(\alpha_i)$ is
increasing to $p(\alpha )$, we are done.

Now we may assume that $p(\alpha_i)=\omega^v\cdot k$ for all $i$ and
$p(\alpha )=\omega^v\cdot (k+1)$.  (This argument will also apply to
the successor
ordinal case.)  We have that for each $s\in \Cal O^{\omega^v\cdot
k}\left(\delta,\left(f_n\right),K\right)$ there is
an infinite set $L_s$ and ordinals $\gamma_s$, $\lambda_s$ such that
$$s\in \Cal O_{+}^{\gamma_s}\left(\delta,\left(f_n\right)_{n\in L_s}
,K\right)\cap \Cal O_{-}^{\lambda_s}\left(\delta,\left(f_n\right)_{
n\in L_s},K\right)$$
and
$$\text{$\min\left(\gamma_s+\lambda_s,\lambda_s+\gamma_s\right)\ge
\omega^v\cdot k$}.$$
Let $S=\left\{s(m):m\in \Bbb N\right\}$ be a countable dense subset of
$
\Cal O^{\omega^v\cdot k}\left(\delta,\left(f_n\right),K\right)$.
 For each $m\in \Bbb N$ there is an infinite set $L_m\subset L_{m-
1}\subset \Bbb N$ and ordinals
$\gamma_m$ and $\lambda_m$ such that
$$s(m)\in \Cal O_{+}^{\gamma_m}\left(\delta,\left(f_n\right)_{n\in
L_m},K\right)\cap \Cal O_{-}^{\lambda_m}\left(\delta,\left(f_n\right
)_{n\in L_m},K\right)$$
and
$$\min\left(\gamma_m+\lambda_m,\lambda_m+\gamma_m\right)\ge\omega^
v\cdot k.$$
By a diagonalization argument we may assume that the set $L$ does
not depend on $m$ and that $\gamma_m=\omega^v\cdot k_m$ and $\lambda_
m=\omega^v\cdot j_m$.  Hence
$$\Cal O^{\omega^v\cdot k}\left(\delta,\left(f_n\right),K\right)\subset
\cup_{i+j=k}\Cal O_{+}^{\omega^v\cdot i}\left(\delta,\left(f_n\right
)_{n\in L},K\right)\cap \Cal O_{-}^{\omega^v\cdot j}\left(\delta,\left
(f_n\right)_{n\in L},K\right).$$
Because each of these sets is closed, we need only consider those
sets which contain $t.$ Moreover observe that if
$$t\in \Cal O_{+}^{\omega^v\cdot i}\left(\delta,\left(f_n\right)_{n
\in L},K\right)\cap \Cal O_{-}^{\omega^v\cdot j}\left(\delta,\left(
f_n\right)_{n\in L},K\right)$$
for two different pairs $(i,j)$ and $(i^{\prime},j^{\prime})$ then
$$t\in \Cal O_{+}^{\omega^v\cdot i^{\prime\prime}}\left(\delta,\left
(f_n\right)_{n\in L},K\right)\cap \Cal O_{-}^{\omega^v\cdot j^{\prime
\prime}}\left(\delta,\left(f_n\right)_{n\in L},K\right)$$
where $i^{\prime\prime}=\max\text{ }(i,i^{\prime})$ and $j^{\prime
\prime}=\max\text{ }(j,j^{\prime})$ and $\omega^v\cdot i^{\prime\prime}
+\omega^v\cdot j^{\prime\prime}\ge p(\alpha )$,
as required.  If there is only one such pair $(i,j)$, we may assume
that a neighborhood of $t$ (relative to $\Cal O^{\omega^v\cdot k}\left
(\delta,\left(f_n\right),K\right)\big)$ is contained in
$\Cal O_{+}^{\omega^v\cdot
i^{\prime\prime}}\left(\delta,\left(f_n\right
)_{n\in L},K\right)\cap \Cal O_{-}^{\omega^v\cdot
j^{\prime\prime}}\left
(\delta,\left(f_n\right)_{n\in L},K\right)$.  We have that
$$t\in \Cal O^{\omega^v}\left(\delta,\left(f_n\right)\right),\Cal O_{
+}^{\omega^v\cdot i^{\prime\prime}}\left(\delta,\left(f_n\right)_{n
\in L},K\right)\cap \Cal O_{-}^{\omega^v\cdot j^{\prime\prime}}\left
(\delta,\left(f_n\right)_{n\in L},K\right)\big)$$
and thus the case $\alpha =\omega^v$ gives that
$$t\in \Cal O_{\epsilon}^{\omega^v}\left(\delta,\left(f_n\right)_{n
\in L^{\prime}},\Cal O_{+}^{\omega^v\cdot
i^{\prime\prime}}\left(\delta,\left
(f_n\right)_{n\in L},K\right)\cap \Cal O_{-}^{\omega^v\cdot j^{\prime
\prime}}\left(\delta,\left(f_n\right)_{n\in L},K\right)\right)$$
for some infinite $L^{\prime}\subset L$ and $\epsilon =+$ or $-$.
Hence

$$t\in \Cal O_{+}^{\omega^v\cdot i}\left(\delta,\left(f_n\right)_{n
\in L^{\prime}},K\right)\cap \Cal O_{-}^{\omega^v\cdot
j}\left(\delta,\left
(f_n\right)_{n\in L^{\prime}},K\right),$$
where $i=i^{\prime\prime}+1$ and $j=j^{\prime\prime}$ if $\epsilon
=+$, and $i=i^{\prime\prime}$ and $j=j^{\prime\prime}$ $+1$ if
$\epsilon
=-$,
as required.\qed\enddemo

Now we are ready to begin the proof of Theorem 3.1.

\demo{Proof}  By Lemma 3.4 we need only prove it in the simpler
case indicated above, i.e., suppose that $\Cal O_{+}(\epsilon )\ge
\alpha$.  Let $f(t)=\lim f_n(t)$.
Fix $\rho ,$ $\epsilon /8>\rho >0$.  For any $g\in C(K)$ let $C(g,
k_0)=\left\{k:|g(k)-f(k_0)|<\rho\right\}$.

The proof is by induction on $\alpha$.  As usual we must actually prove
a
little stronger statement to make the induction work.\newline
INDUCTIVE HYPOTHESIS:  Suppose that $k_1,k_2,\ldots ,k_j$ are a finite
number
of points in $\Cal O_{+}^{\alpha}\left(\epsilon ,\left(f_
n\right),K\right)$ and $\Cal N_1,\Cal N_2,\ldots ,\Cal N_j$ are
neighborhoods of
$k_1,k_2,\ldots ,k_j$, respectively, then for every $\beta <\alpha$
there is a tree $
\Cal T$ on
$\left(f_n-f_m\right)$ of order $\beta$ such that for $i=1,2,\ldots
,j$,
$$\Cal T_{\Cal N_i}=\left\{\left(x_{1|_{\Cal N_i}},x_{2|_{\Cal N_i}}
,\ldots ,x_{n|_{\Cal N_i}}\right):\left(x_1,x_2,\ldots ,x_n\right)
\in \Cal T\right\},$$
is an $\ell^1$ tree of order $\beta$ as in Lemma 3.3, i.e., for each $
i$ and branch
there is a mapping into the pairs of subsets of $K$ satisfying a) and
b)
with $\delta =\epsilon /2$ and $r=\epsilon /4$.  Moreover we may assume
that for each
pair of sets
$$\rho_{\Cal B}\left(f_{n_1}-f_{m_1},\ldots ,f_{n_j}-f_{m_j}\right
)=\left(A_j,B_j\right),$$
$$A_j\subset\cap_{i=1}^{j-1}C\left(f_{n_i},k\right)\cap C\left(f_{
m_i},k\right)\cap C\left(f_{m_j},k\right)\cap\left\{s:|f_{n_j}(s)-
f_{n_j}(k)|<\rho\right\},$$
and
$$B_j\subset\cap^j_{i=1}C\left(f_{n_i},k^{\prime}\right)\cap C\left
(f_{m_i},k^{\prime}\right)\text{ for some $k$ and $k^{\prime}$ in $
K$.}$$

Assume the inductive hypothesis holds for all $\beta <\alpha$.  First
suppose
that $\alpha$ is not a limit ordinal and that $\Cal O_{+}^{\alpha}$
$\left
(\epsilon ,\left(f_n\right),K\right)\neq\emptyset$.  Let $\left(k_
i\right)$ and
$\left(\Cal N_i\right)$ be as above.  For each $i$ there is an $N_
i\in \Bbb N$ such that for each
$n\ge N_i$, there is an $M_n^i$ such that
$$\cap_{m\ge M_n^i}A_{n,m}^{+}\cap \Cal O_{+}^{\alpha -1}\left(\epsilon
,\left(f_n\right),K\right)\cap \Cal N_i\neq\emptyset .$$
Let $n\ge\max\left\{N_i\right\}$ such that $|f_n\left(k_i\right)-f\left
(k_i\right)|<\rho$, let $M=\max\left\{M^i_n\right\}$ and
choose a point $\displaystyle{k^{\prime}_i\in\cap_{m\ge M}A_{n,m}^{
+}\cap \Cal O_{+}^{\alpha -1}\left(\epsilon ,\left(f_n\right),K\right
)\cap \Cal N_i}$, for
each $i$.  We may also assume that for all $p\ge M$,
$$|f_p\left(k_i\right)-f\left(k_i\right)|<\rho$$
and
$$|f_p\left(k^{\prime}_i\right)-f\left(k^{\prime}_i\right)|<\rho ,$$
for $i=1,2,\ldots ,j$.  Fix $m\ge M$ and for each $i$ let $\Cal N^{
\prime}_i$ be a neighborhood of
$k^{\prime}_i$ contained in
$$\Cal N_i\cap A_{n,m}^{+}\cap C\left(f_m,k^{\prime}_i\right)\cap\left
\{k:|f_n(k)-f_n\left(k^{\prime}_i\right)|<\epsilon /8\right\}$$
and let $\Cal N^{\prime\prime}_i$ be a neighborhood of $k_i$ contained
in
$$\Cal N_i\cap C\left(f_m,k_i\right)\cap C\left(f_n,k_i\right).$$

Now by the inductive hypothesis for every $\beta <\alpha$ there is a
tree $
\Cal T$
on $\left(f_p-f_q\right)$ such that $\Cal T_{\Cal N^{\prime}_i}$ and $
\Cal T_{\Cal N^{\prime\prime}_i}$ are $\ell^1$ trees of order $\beta$
for
$i=1,2,\ldots ,j$, and satisfy a) and b) of the lemma with $\delta
=\epsilon /2$ and
$r=\epsilon /4$.  Let
$$\Cal T^{\prime}=\left\{\left(f_n-f_m,x_1,x_2,\ldots ,x_n\right):\left
(x_1,x_2,\ldots ,x_n\right)\in \Cal T\right\}.$$
Clearly this is a $\left(\beta +1\right)$ tree on $\left(f_n\right)$.
We need to check the hypothesis
of the lemma for $\Cal T^{\prime}_{\Cal N_i}$, $i=1,2,\ldots ,j$.  To
define $
\rho^{\prime}_{\Cal B^{\prime}}$ on a branch
$$\Cal B^{\prime}=\left\{\left(f_n-f_m\right),\left(f_n-f_m,x_1\right
),\left(f_n-f_m,x_1,x_2\right)\ldots\right\}\text{ of }\Cal
T^{\prime}_{
\Cal N_i},$$
we let
$$\rho^{\prime}_{\Cal B^{\prime}}\left(f_n-f_m,x_1,x_2,\ldots
,x_n\right
)=\rho_{\Cal B}\left(\left(x_1,x_2,\ldots ,x_n\right)\right)$$
and
$$\rho^{\prime}_{\Cal B^{\prime}}\left(\left(f_n-f_m\right)\right)
=\left(\Cal N^{\prime}_i,\Cal N^{\prime\prime}_i\right),$$
where $\rho_{\Cal B}$ denotes the mapping from the branch $\Cal B=\left
\{\left(x_1\right),\left(x_1,x_2\right),\ldots\right\}$
of $\Cal T_{\Cal N_i}$.  If $r=\epsilon /4$ and $\delta =\epsilon
/2$, the hypothesis of the lemma is
satisfied and thus $\Cal T^{\prime}_{\Cal N_i}$ is an $\ell^1$-tree of
order $
\beta +1$.  Because this is
true for all $\beta <\alpha$ we get an $(\alpha +1)$ $\ell^1$-tree.

For limit ordinals the conclusion is obvious.  (Note we can actually
get an $(\alpha +1)$ $\ell^1$-tree in this case as well.)

The moreover assertion allows us to conclude that the we can
construct a $\ell^1$-tree on $\left(f_n\right)$ of the same order.
Indeed we claim that
the tree obtained by replacing in each coordinate $f_n-f_m$ by $f_
n$ is the
required tree.  First if $f$ is continuous we can choose a neighborhood
$\Cal N$ of the point $k_0\in \Cal O^{\alpha}\left(\epsilon ,\left
(f_n\right),K\right)$ so that $\Cal
N\subset\left\{k:|f(k)-c|<\rho\right
\}$, where
$c=f\left(k_0\right)$, and restrict all of the functions to this set $
\Cal N$.  The proof
then shows that the sets $\left\{k:f_{n_i}\ge c+\epsilon -\rho\right
\}$ and $\left\{k:f_{n_i}\le c+\rho\right\}$ are
Boolean independent for any $\left(f_{n_1}-f_{m_1},f_{n_2}-f_{m_2}
,\ldots ,f_{n_j}-f_{m_j}\right)$ in the
tree constructed.

If $f$ is not continuous we use the following lemma.

\proclaim{Lemma 3.5}  If $\left(x_n\right)$ is a uniformly bounded
sequence in a
Banach space $X$ such that $||\sum a_nx_n||\ge\delta\sum |a_n|$ for all
$\left
(a_n\right)\in \Bbb R^{\Bbb N}$ and
$y\in X$, $||y||\le 1$, then
$$||\sum a_n\left(x_n+y\right)||\ge\delta^{\prime}\sum |a_n|,$$
where $\delta^{\prime}=\max\left\{\delta
d\left(y,\left[x_n\right]\right
)/2,\delta -||y||\right\}$.\endproclaim

We omit the simple proof.

Our tree was actually constructed using sequences $\left(f_{n_i}-f_{
m_i}\right)_{i=1}^N$ and
sets\linebreak $\left(A_{n_i},B_{n_i}\right)$ where
$|f_{m_j}(k)-f(k^{\prime}
)|<\rho$ for all
$\displaystyle{k\in\cap^N_{i=1}A^{\epsilon_i}_{n_i}}$, for all $j$, for
some\linebreak
$\displaystyle{\text{$k\left(\epsilon_1,\epsilon_2,\ldots ,\epsilon_
N\right)\in$ \linebreak$\cap^N_{i=1}A^{\epsilon_i}_{n_i}$}}$.  Thus we
can replace each
$f_{m_i}$ by $g$ where
$$g (k)=\cases f\left(k\left(\epsilon_1,\epsilon_2,\ldots ,\epsilon_
N\right)\right)&\text{if }k\in\cap^N_{i=1}A^{\epsilon_i}_{n_i},\\
0&\text{otherwise}\endcases$$
with a loss of at most $\rho$.  Observe that because $\lim f_n=f$,
$\displaystyle{\lim_N\inf d\left(f,\left[f_n-f:n\ge N\right]\right
)}$\linebreak $\ge||f||$.  Thus we get the
estimate
$$\delta^{\prime}\ge\max\left\{\left(\frac {\epsilon}2-\rho\right)
||f||,\epsilon -||f||-\rho\right\}$$
from the lemma.\qed\enddemo

Let us now examine the relationship between the oscillation index and
$\ell^1$-spreading models.

\proclaim{Proposition 3.6}  Let $\left(f_n\right)$ be a pointwise
converging
sequence on a Polish space $(K,d)$ with limit 0.  If $k\in \Cal O^{
\omega}\left(\epsilon ,\left(f_n\right),K\right)$
and $\Cal N$ is a neighborhood of $k$, then for every infinite
$M\subset
\Bbb N$ there is
an $L\subset M$ such that $\left(f_{n|\Cal N}\right)_{n\in L}$ has $
\ell^1$ spreading model constant at \linebreak least
$\epsilon$.\endproclaim

\demo{Proof}  Fix $j\in \Bbb N$ and choose $L\subset M$ such that
$\left
(f_{n|\Cal N}\right)_{n\in L}$ has a
spreading model.  Because $k\in \Cal O^{2j}\left(\left(f_n\right),
K,\delta\right)$, a careful examination of
the proof of the previous theorem shows that there is an $\epsilon
=+$ or $-$
such that given any $N$ there is an $\Cal
L\subset\left\{N,N+1,\ldots\right
\}$ with cardinality $j$
with $\Cal N\cap\cap_{\ell\in \Cal L}A_{\ell
,m}^{\epsilon}\neq\emptyset$, for all large $
m$.  Therefore
$\left(f_{n|\Cal N}\right)_{n\in L}$ has $\ell^1$ spreading model
constant at least $
\epsilon$.\qed\enddemo

We wish to present two more examples, but before doing so we will
introduce another ordinal index which we call the spreading model
index and show that it is closely related to the oscillation index.
This index is defined in terms of the spreading model constant and
should be compared with the averaging index.

Let $S^0\left(\epsilon ,\left(f_n\right),K\right)=K$ and assuming that
$S^{\alpha}\left(\epsilon ,\left(f_n\right),K\right)$ has been
defined let
$$\align S^{\alpha +1}\left(\varepsilon
,\left(f_n\right),K\right)=&\big\{
x^{*}:\text{ for every neighborhood }\Cal N\text{ of }x^{*}\text{
relative to}\\
&S^{\alpha}\left(\varepsilon ,\left(f_n\right),K\right)\text{ and
infinite }
M\subset \Bbb N\text{ there is an }L\subset N\\
&\text{with }\ell^1-SP\left(\left(f_{n|\Cal N}\right)\right)_{n\in
L}\ge\varepsilon\big\}.\endalign$$

As usual if $\alpha$ is a limit ordinal
$\displaystyle{S^{\alpha}\left(\epsilon ,\left(f_{_{}n}\right),K\right
)=\cap_{\beta <\alpha}S^{\beta}\left(\epsilon ,\left(f_n\right),K\right
)}$ and the ordinal
index will be denoted by $o(S,\epsilon )$.

Next we will show that the oscillation index and the spreading model
index measure essentially the same thing.

\proclaim{Proposition 3.7}  Suppose that $\left(f_n\right)$ is a weakly
null
sequence in $B_{C(K)}$ for some compact metric space $K$.  Then

\itemitem{i)}  if $S^1\left(\epsilon ,\left(f_n\right),K\right)\neq
\emptyset ,\Cal O^1\left(\epsilon^{\prime},\left(f_n\right),K\right
)\neq\emptyset$ for every $\epsilon^{\prime}<\epsilon$.
\itemitem{ii)} if $\Cal O^{\omega}\left(\epsilon ,\left(f_n\right)
,K\right)\neq\emptyset$, for some $\epsilon >0$, $S^1\left(\epsilon
,\left(f_n\right),K\right)\neq\emptyset$.
\endproclaim

\demo{Proof}  Suppose that $\Cal O^1\left(\epsilon^{\prime},\left(
f_n\right),K\right)=\emptyset$.  Then for each $k\in K$
there is a neighborhood $\Cal N_k$ of $k$ and a subsequence $\left
(f_n\right)_{n\in M_k}$ such
that $||f_{n|\Cal N_k}||\le\epsilon^{\prime}$  for all $n\in M_k$.
Clearly
$\displaystyle{||\sum_{n\in F}f_{n|\Cal N_k}||\le
|F|\epsilon^{\prime}}$ for every finite subset $
F$ of $M_k$.
Hence $S^1\left(\epsilon ,\left(f_n\right),K\right)=\emptyset$ for
every $
\epsilon >\epsilon^{\prime}$.

The second assertion follows from Proposition 3.6.\qed\enddemo

\proclaim{Corollary 3.8}  Suppose that $\left(f_n\right)$
is a weakly null sequence
on a compact metric space $K$.  Then

\itemitem{i)} If $o(S,\epsilon )\ge\alpha$, then $\Cal O(\epsilon^{
\prime})\ge\alpha$, for every $\epsilon^{\prime}<\epsilon$.
\itemitem{ii)} If $\Cal O(\epsilon )\ge\omega^{1+\alpha}$, then $o
(S,\epsilon )\ge\omega^{\alpha}$.
\endproclaim

\example{Example 4}  $L^1$.  Of course $\ell\left(L^1,1\right)=\omega_
1$.  However if $\left(f_n\right)$ is
a weakly convergent sequence in $L^1$,
then $\left(f_n\right)$ is uniformly integrable
and hence by Dor's Theorem [Dor] for every $K<\infty$ there is an $
n$ such
that if $F$ is a set of integers of cardinality $n$, $\left(f_j\right
)_{j\in f}$ is not $K$
equivalent to the unit vector basis of $\ell_n^1$.  Therefore
$\ell^1$ $-SP\left(f_n\right)_{n\in L}=0$ for all $L$ and thus $\Cal
O^{
\omega}\left(\epsilon ,\left(f_n\right),K\right)=\emptyset$ for every
$\epsilon >0$.\endexample

\example{Example 5}  The spaces Szlenk [Sz] used to show that there
are reflexive spaces with arbitrarily large (countable) Szlenk index
are defined inductively as $X_1=\ell^2$, $X_{\alpha
+1}=\left(X_{\alpha}
\oplus\ell^2\right)_1$, and for a limit
ordinal $\alpha$, $\displaystyle{X_{\alpha}-\big(\sum_{\beta <\alpha}
X_{\beta}\big)_2}$.  The $\ell^1$-index of $X_{\alpha}=\alpha$
for $\epsilon =1$ and increases to $\alpha\omega$
as $\epsilon$ goes to 0.\endexample

Now let us consider the oscillation index.  First suppose that $\alpha$
is a
limit ordinal.  If $\left(f_n\right)$ is a weakly null sequence in $
X_{\alpha}$, then by
passing to a subsequence (which can only increase the index) we may
assume that $f_n=g_n+h_n$ where $\displaystyle{g_n\in\sum_{\beta\le
\lambda}X_{\beta}}$ for all $n$
for some $\lambda <\alpha$ and $\displaystyle{h_n\in\sum_{\beta\in
B_n}X_{\beta}}$ where
$B_n=\left\{\beta :\beta_n<\beta\le\beta_{n+1}\right\}$ for all $n$ and
$
\beta_1>\lambda$.  Note that $\left[h_n\right]=\ell^2$ and
thus the $\ell^1$ spreading model index of $X_{\alpha}$ is the supremum
of the
indices of $X_{\beta}$, $\beta <\alpha$.  Now suppose that $\alpha
=\beta +k$ for some integer $n$.
Then $X_{\alpha}=\left(X_{\beta}\oplus\undersetbrace k\to{\ell^2\oplus
\ell^2\oplus\cdots\oplus\ell^2}\right)_1$, and we may write $f_n=g_
n+h_n$ where
$g_n\in X_{\beta}$ and $\displaystyle{h_n\in\sum^k_1\oplus\ell^2}$ for
all $
n$.  However we again
have that $\ell^1-SP\left(f_n\right)=\ell^1-SP\left(g_n\right)$
and hence the spreading model
index of $\left(f_n\right)$ is the spreading model index of $\left
(g_n\right)$.  Thus the
spreading model index of $X_{\alpha}$ is the same as the index of $
\ell^2$, for all
$\alpha$, namely 0.
\vskip .15in
This last family of examples illustrates the fact that the oscillation
index and the spreading model index are really measuring something
stronger than the existence of many $\ell_n^{1\prime}s$ in a space.

Next we will examine the ordinal index of a Baire-1 function $f$ on a
Polish space $(K,d)$.  The index is defined by considering two real
numbers $c$ and $d$, $c<d$, and the disjoint $G_{\delta}$ sets,
$$C=\left\{k\in K:f(k)\le c\right\}\text{ and }D=\left\{k\in K:f(k
)\ge d\right\}.$$
$L(f,c,d)$ is the smallest ordinal $\alpha$ such that there is a
decreasing
family of closed sets $F_{\beta},$ $\beta\le\alpha$, with $F_0=K$, $
F_{\alpha}=\emptyset$, and for all
$\beta <\alpha$, $F_{\beta}\backslash F_{\beta +1}$ is disjoint from $
C$ or from $D$ and at a limit ordinal $\gamma$,
$F_{\gamma}=\cap_{\beta <\gamma}F_{\beta}$.
(See [Bo] where the definition is given
in complementary terms or [K,p452].)

Bourgain shows that if $\left(f_n\right)$ is a pointwise converging
sequence of
continuous functions with limit $f$ and $\epsilon <(d-c)/2$ then $
\omega^{\ell\left(\left[f_n\right],\epsilon\right)+1}$
is greater than $L(f,c,d)$.  (Actually his result gives a slightly
smaller
bound.)  We wish to show that in fact the oscillation index is also
large.  The proof of the following proposition is similar to that of
Lemma 5, [H-O-R].

\proclaim{Proposition 3.9}  Suppose that $\left(f_n\right)$ is a
pointwise
converging sequence of continuous function on a compact metric
space $K$.  Then if $d-c>\epsilon$ and if $L(f,c,d)=\beta +m$, where $
\beta$ is a limit
ordinal and $m<\omega$, then $\Cal O_{+}\left(\epsilon ,\left(f_n\right
),K\right)\ge\beta +(m-1)/2$.\endproclaim

\demo{Proof}  Consider the following family of closed sets where
$C=\left\{k:f(k)\le c\right\}$ and $D=\left\{k:f(k)\ge d\right\}$, $
F_0=K,$ $F_1=\overline D$, $F_2=\overline {F_1\cap C}$,
$F_3=\overline {F_2\cap D}$, and in general, $F_{\alpha
+2n+1}=\overline {
F_{\alpha +2n}\cap D}$ and
$F_{\alpha +2n+2}=\overline {F_{\alpha +2n+1}\cap C}$ if $\alpha$ is
even and $
n\in \Bbb N$.  (Limit ordinals are
even.)  If $\alpha$ is a limit ordinal $F_{\alpha}=\cap_{\beta <\alpha}
F_{\beta}$.

It is easy to see that if $\alpha$ is even $\left(F_{\alpha}\backslash
F_{\alpha +1}\right)\cap D=\emptyset$ and
$\left(F_{\alpha +1}\backslash F_{\alpha +2}\right)\cap C=\emptyset$.
 Next we will show that
$F_{\alpha +2}\cap D\subset \Cal O_{+}^1$ $\left(\epsilon ,\left(f_
n\right),F_{\alpha +1}\right)$ for all $\alpha$ even.

Let $d^{\prime}\in F_{\alpha +2}\cap D$.  This implies that $d^{\prime}
\in\overline {F_{\alpha +1}\cap C}\cap D$.  $\overline {F_{\alpha
+1}\cap C}\cap D$
contains points from $D$ which are (non-trivial) limits of points in
$F_{\alpha +1}\cap C$ and hence there exists a sequence
$\left(c_k\right
)$ in $F_{\alpha +1}\cap C$ with
limit $d^{\prime}$.  Moreover we may assume that no $c_k$ is in $\Cal
O_{
+}^1\left(\epsilon ,\left(f_n\right),F_{\alpha}\right)$
(else $d^{\prime}$ would be also).  Choose $N\in \Bbb N$ such that $
|f_n(d^{\prime})-f(d^{\prime})|$ is less
than $\delta /4$ for all $n\ge N$, where $0<\delta <(d-c)-\epsilon$.
For each $
k\in \Bbb N$ there
is an $M_k$ such that $f_m(c_k)<\delta /4$ $+f(c_k)$ for all $m\ge
M_k$.  Now if $\Cal N$ is a
neighborhood of $d^{\prime}$ and $n\ge N$ then for some $L\in \Bbb N$,
$
c_k\in \Cal N$ and
$f_n(c_k)>f(d^{\prime})-\delta /4$, for all $k\ge L$.  Because
$$f_n(c_k)-f_m(c_k)>f(d^{\prime})-\delta /4-f(c_k)-\delta /4\ge d-
c-\delta /2,$$
for all $k\ge L$ and $m\ge M_k,$ $c_k\in \Cal N\cap\cap_{m\ge M_k}
A_{n,m}^{+}$.  Hence
$d^{\prime}\in \Cal O_{+}^1\left(\epsilon ,\left(f_n\right),F_{\alpha
+1}\right)$ and $F_{\alpha +3}=\overline {F_{\alpha +2}\cap D}\subset
\Cal O_{+}^1\left(\epsilon ,\left(f_n\right),F_{\alpha +1}\right)$, for
all
$\alpha$ even.

Because $\Cal O_{+}^{\alpha +1}\left(\epsilon ,\left(f_n\right),K\right
)=\Cal O_{+}^1\left(\epsilon ,\left(f_n\right),\Cal O^{\alpha}_{+}
(\varepsilon ,(f_n),K\big)\right)$.  A simple induction
argument\linebreak shows that $F_{\alpha +2m+1}\subset \Cal O_{+}^
m\left(\epsilon ,\left(f_n\right),F_{\alpha +1}\right)$ for any integer
$
m$
and even ordinal $\alpha$.  It follows then that if $\beta$ is a limit
ordinal
$F_{\beta}\subset \Cal O_{+}^{\beta}\left(\epsilon ,\left(f_n\right
),K\right)$ and for any $m\in \Bbb N$, $F_{\beta +2m+1}\subset \Cal O_{
+}^{\beta +m}\left(\epsilon ,\left(f_n\right),K\right)$.\qed\enddemo

\proclaim{Corollary 3.10}  Suppose that $f$ is a Baire-1 function on a
compact metric space $K$ and that $\left(f_n\right)$ is
a sequence of continuous
functions on $K$ with $||f_n||\ge 1$ which converge to $f$ pointwise.
If
$L(f,c,d)=\beta +m$, where $c<d$, $\beta$ is a limit ordinal and $
m\in \Bbb N$, then for
any $\varepsilon <d-c$ there is an $\epsilon /2\text{ }\ell^1$-tree of
order $
\beta +m/2$ on $\left(f_n-f_m\right)$, and there is an $\epsilon
/2\text{ }
\ell^1$-tree on
$(f_n)$ of the same order.\endproclaim

\demo{Proof}  From the above proposition we get that
$\Cal O_{+}^{\beta +(m-1)/2}\left(\epsilon ,\left(f_n\right),K\right
)\neq\emptyset$.  The proof of Theorem 3.1 shows that
there is an $\ell^1$-tree on $\left(f_n-f_m\right)$ of order $\beta
+(m-1)/2+1$ with lower
estimate $\epsilon$.  The second assertion follows from an examination
of the
proofs of Proposition 3.9 and Theorem 3.1.  It is easy to see that in
the proof of the theorem we can replace $A^{+}_{n,m}$ by $\left\{k
:f_n(k)-c>\epsilon +\rho\right\}$,
where $\rho <d-c-\epsilon$, and always choose the points $k_i$ and $
k^{\prime}_i$ from $C$.
We then get that the sets $\left(\left\{k:f_{n_s}(k)>c+\epsilon
+\rho\right
\},\left\{k:f_{n_s}(k)<c+\rho\right\}\right)_{s=1}^j$
are Boolean independent for any node $\left(f_{n_1},f_{n_2},\ldots
,f_{n_j}\right)$ of the tree
constructed.\qed\enddemo

\head 4.  Construction of weakly null sequences with large oscillation
index\endhead

In this section we wish to generalize the construction of Schreier of
a weakly null sequence $\left(x_n\right)$ with no subsequence having
the
Banach-Saks property.  As observed by Pelczynski and Szlenk [P-S]
the  Schreier sequence is a 1-suppression unconditional basis in
$C(\omega^{\omega})$.  Our goal is to prove

\proclaim{Theorem 4.1}  For every $\alpha <\omega_1$ there is a weakly
null
sequence $\left(x^{\alpha}_n\right)$ in
$C\left(\omega^{\omega^{\alpha}}\right
)$ with \linebreak$\Cal O^{\omega^{\alpha}}\left(1-\epsilon ,\left
(x^{\alpha}_n\right),\omega^{\omega^{\alpha}}\right)\neq\emptyset$ for
every
$\epsilon >0$.  Moreover for each $\alpha$, $\left(x^{\alpha}_n\right
)$ is a sequence of indicator
functions and $\left(x^{\alpha}_n\right)$ is a
1-suppression unconditional basic sequence.\endproclaim

To construct these sequences and verify their properties it is useful
to have several different viewpoints.  The first viewpoint is
contained in the following result.

\proclaim{Proposition 4.2}  Let $\Cal F$ be a family of finite subsets
of $
\Bbb N$
such that
\itemitem{i)}  if $F\in \Cal F$ and $G\subset F$
 then $G\in \Cal F$, i.e., the family is
adequate.
\itemitem{ii)}  $\left\{n\right\}\in \Cal F$ for all $n\in \Bbb N$
\itemitem{iii)}  if $F_j\in \Cal F$ for $j=1,2,\ldots$ and $1_{F_j}$
converges pointwise
to $1_F$ then $F\in \Cal F$.
\item{} Let $x_n=1_{\left\{F\in \Cal F:n\in F\right\}}$ for
$n=1,2,\ldots$ .
\item{} Then $\Cal F$ is a countable compact metric space under the
topology
induced by identifying $\Cal F$ with $\left\{1_F:F\in \Cal F\right
\}$ under the topology of
pointwise convergence.  $\left(x_{_{}n}\right)$ is a weakly null
1-suppression
unconditional basic sequence in $C(\Cal F)$.
\endproclaim

\demo{Proof}  The first assertion is an easy consequence of iii) and
we omit the argument.  Note that $x_n(F)\neq 0$ if and only if $n\in
F$.
Hence
$$||\sum a_nx_n||=\sup\left\{|\sum_{_{n\in F}}a_n|:F\in \Cal F\right
\}.$$
If $G\subset \Bbb N$, then by i)
$$||\sum_{n\in G}a_nx_n||=\sup\left\{|\sum_{n\in G\cap F}a_n|:F\in\!
\Cal F\right\}=\sup\left\{|\sum_{n\in F}a_n|:F\in \Cal F\text{ and }
F\subset G\right\}$$
Clearly this is not greater than $||\sum a_nx_n||$.  Because each $
F\in \Cal F$ is
finite, $x_n(F)\neq 0$ for only finitely many $n$.  Thus $\left(x_
n\right)$ is weakly null.\qed\enddemo

This proposition gives us an easy way of defining and verifying the
properties of the Schreier sequence.  As in the previous section let
$$\Cal F_1=\left\{F\subset \Bbb N:\min F\ge\text{ card }F\right\}$$
(We consider $\emptyset$ to be in $\Cal F_1$.)  If $(x_n)$ is defined
as in the proposition
then it follows that $(x_n)$ is a 1-unconditional basic sequence and is
weakly null.  Finally if $L\subset \Bbb N$ is infinite and for each $
k\in \Bbb N$, we let
$L_k$ be the first $k$ elements of $L$
$$||\sum_{n\in L_{2k}}x_n||\ge k$$
because $L_{2k}\backslash L_k\in \Cal F_1$ and
$$||\sum_{n\in L_{2k+1}}x_n||\ge k$$
because $L_{2k+1}\backslash L_{k+1}\in \Cal F_1$.  Hence
$$||\sum_{n\in L_k}x_n||/k\ge (k-1)/(2k)\text{ for all }k,$$
and therefore $\left(x_n\right)_{n\in L}$ fails the Banach-Saks
property.  (See [D1,p.78].)

The major drawback to this representation of the Schreier sequence
is that it is difficult to understand the topology of $\Cal F_1$ and
its
relation to $\left(x_n\right)$.  In this section we will
prove some results twice.
First we will give proofs based on representations like that above
for the Schreier sequence.  The second will be given using trees.  We
have found that this second viewpoint is more intuitive (It is easy to
draw pictures of the trees.) and we used it to establish these results
originally.

Our next goal then is to use trees to describe the Schreier sequence
and the underlying topological space and in particular to show that
this sequence and its generalizations could be obtained by beginning
with the coordinate functions on the Cantor set, and then essentially
restricting them to suitable subsets of the Cantor set.  More
precisely, if we let $C=\left\{-1,1\right\}^{\Bbb N}$ and $r_n\left
(\left(\epsilon_i\right)\right)=\epsilon_n$, then
$x_n=\left(r_n+1\right)/2$ is a sequence of indicator functions.  If $
K$ is a
compact subset of $C$, then the sequence $(x_{n|K})$ is a sequence of
indicator functions in $C(K)$ which will be equivalent to the Schreier
sequence if $K$ is properly chosen.

Now let us work backwards from a sequence of indicator functions
to find an appropriate minimal underlying topological space.  Suppose
that $(x_n)$ is a sequence of indicator functions on a set $K$.  Define
a
tree $\Cal T$ on $\{-1,1\}$ by
$$\Cal T=\cup^{\infty}_{n=1}\left\{\left(\epsilon_1,\epsilon_2,\ldots
,\epsilon_n\right):\cap^n_{i=1}(\text{supp }x_i)^{\epsilon_i}
\neq\emptyset\right
\}$$
where
$$\text{supp }x_i=(\text{supp }x_i)^1=\left\{k\in K:x_i(k)=1\right
\}$$
and
$$(\text{supp }x_i)^{-1}=K\backslash (\text{supp }x_i)^1.$$
Let $\overline {\Cal T}=\Cal T\cup\left\{\left(\epsilon_1,\epsilon_
2,\ldots\right):\cap^n_{i=1}(\text{supp }x_i)^{\epsilon_i}\neq
\emptyset\text{ for all }
n\right\}$.  $\overline {\Cal T}$ is also a tree
and both are subtrees of the full dyadic tree
$$\Cal D=\cup^{\infty}_{n=1}\{-1,1\}^n\cup \{-1,1\}^{\Bbb N}.$$
$\Cal D$ has a natural topology given by coordinate-wise convergence.
Of
course in this topology $\{-1,1\}^{\Bbb N}$ is just the
Cantor set and each finite
sequence is an isolated point.  Also any infinite sequence is the limit
of its restrictions to the first $n$ coordinates, i.e., the nodes above
it.
 Note also that this tree is closely related to the Boolean
independence tree $\Cal T\left\{\left(\left(\text{supp }x_n\right)^
1,\left(\text{supp }x_n\right)^{-1}\right)\right\}$.

\proclaim{Lemma 4.3}  $\overline {\Cal T}$ is the closure of $\Cal T$
in $
\Cal D$.
\endproclaim

\demo{Proof}  Obvious.\qed\enddemo

{\it Throughout the remainder of this section we will assume that the
sequence $(x_n)$ is pointwise convergent to 0 on K, K is a compact
metric space and each $x_n$ is continuous.\/}

Given such a sequence $(x_n)$ the tree $\Cal T((x_n))$ defined above
will be
called the {\it tree associated to $(x_n)$\/}.  Because $K$ is compact
and $
(x_n)$ is
weakly null, $\overline {\Cal T}$ contains no elements
with infinitely many coordinates
equal to 1 and every node is on an infinite branch.  Hence $\overline {
\Cal T}$ is
countable and therefore homeomorphic to some countable ordinal in
the order topology.  Unfortunately because the tree is not well
founded, one cannot use the order of the tree $\Cal T$ to determine the
topological type.  To get around this problem we must study the
relationship between the topology of $\overline {\Cal T}$ and
the structure of the tree.

Before we embark on the study of the topology let us consider a
property of trees is analogous to property i) of Proposition 4.2.  In
an unpublished paper H.P. Rosenthal introduced the following
notation.

\definition{Definition}  Let us say that a tree $\Cal T$ on $\{-1,
1\}$ is {\it weakly
independent\/} if \linebreak$\left(\epsilon_1,\epsilon_2,\epsilon_
3,\ldots ,\epsilon_j\right)\in \Cal T$ implies that for all $\left
(\epsilon^{\prime}_1,\epsilon^{\prime}_2,\epsilon^{\prime}_3,\ldots
,\epsilon^{\prime}_j\right)$
such that $\epsilon^{\prime}_i=-1$ if $\epsilon_i=-1$ and $\epsilon^{
\prime}_i=1$ or -1 if $\epsilon_i=1$,
$\left(\epsilon^{\prime}_1,\epsilon^{\prime}_2,\epsilon^{\prime}_3
,\ldots ,\epsilon^{\prime}_j\right)\in \Cal T$.  Say that a
sequence of indicator functions is
{\it weakly independent\/} if the associated tree is.
\enddefinition

Rosenthal also proved the following result in a slightly different
form.

\proclaim{Proposition 4.4}  A weakly null sequence of (non-zero)
indicator
functions on a compact metric space $K$ determines a weakly
independent tree if and only if it is a 1-suppression unconditional
basic sequence.\endproclaim

\demo{Proof}  First suppose that $(x_n)$ is a sequence of indicator
functions which is a 1-suppression unconditional basic sequence and
$\Cal T$ is the associated tree.  Then for any sequence of real numbers
$(a_n)$ we have
$$||\sum a_nx_n||=\sup\left\{|\sum a_nx_n(k)|:k\in K\right\}=\sup\left
\{\left|\sum_{n:k\in\text{supp }x_n}a_n\right|:k\in K\right\}.$$
In particular suppose that $F\subset \{1,2,\ldots ,j\}$ with
$$\bigcap_{n\in F}(\text{supp }x_n)^{-1}\cap\bigcap\Sb n\notin
F,\\
n\le j \endSb
\text{supp }x_n\neq\emptyset,$$
i.e., the
node $(\epsilon_1,\epsilon_2,\ldots ,\epsilon_j)$ $\in \Cal T$, where $
\epsilon_i=1$ if $n\notin F$ and $n\le j$, and $\epsilon_i=-1$ if
$n\in F$.  In order to check weak independence it is sufficient to
check
the condition on a lower node.  Hence we may assume that if
$\{1,2,\ldots j\}\supset G\supset F$, then $G\not =\left\{1,2,\ldots j
\right\}$, and we must show that
$$\bigcap\Sb n\notin G\\
n\le j \endSb\text{supp}x_n\cap\bigcap_{n\in G}
(\text{supp}x_n)^{-1}\not =\emptyset
.$$
Because $(x_n)$ is 1-suppression unconditional

$$(j-\text{card }G)(1+\text{card }G)=||(1+\text{card }G)
\sum\Sb n\notin G\\ n\le j\endSb x_n+\sum_{n\notin G}(-1)x_n||$$
and thus there exists

$$k\in\bigcap\Sb n\notin G\\ n\le j \endSb(\text{supp }x_n)
\cap\bigcap_{n\in G}(\text{supp }x_n)^{-1}.$$
This implies that the node $(\epsilon_1,\epsilon_2,\ldots ,\epsilon_
j)\in \Cal T$ where $\epsilon_i=1$ if $i\notin G$ and
$i\le j$ and $\epsilon_i=-1$ if $i\in G$.  Thus $\Cal T$ is weakly
independent.

Conversely suppose that $(x_n)$ is a weakly null sequence of indicator
functions and that the associated tree $\Cal T$ is weakly independent.
Then
for any sequence of real numbers $(a_n)$ and finite subset $F$ of $
\Bbb N$, we
claim that
$$||\sum a_nx_n||=\sup\left\{\left|\sum_{n:k\in\text{supp $x_n$}}a_
n\right|:k\in K\right\}$$
$$\ge\sup\left\{\left|\sum_{n\in F:k\in\text{supp $x_n$}}a_n\right
|:k\in K\right\}=||\sum_{n\in F}a_nx_n||.$$
To see the inequality suppose that $k$ is any point in $K$ and
$H=\left\{n:x_n(k)=1\right\}$.  Because $\Cal T$ is weakly independent
there is a point
$k^{\prime}$ in $K$ such that $\left\{n:x_n(k^{\prime})=1\right\}=
F\cap H$.  Hence each sum on the
righthand side of the inequality also occurs on the left.\qed\enddemo

\proclaim{Corollary 4.5}  Suppose that $(x_n)$ is a weakly null
sequence
of non-zero indicator functions in $C(K)$ for some compact metric
space $K$.  Then the following are equivalent.
\itemitem{i)}  The tree associated to $(x_n)$ is weakly independent.
\itemitem{ii)}  $\Cal F=\left\{F\subset \Bbb
N:F=\left\{n:x_n(k)=1\right
\}\text{ for some }k\in K\right\}$ is
adequate.
\itemitem{iii)}  $(x_n)$ is a 1-suppression unconditional basic
sequence.
\endproclaim

\demo{Proof}  We have already shown that $ii)\implies iii)$ and $i
)\Leftrightarrow iii)$.
$i)\implies ii)$ is immediate from the definition of
weakly independent.\qed\enddemo

\remark{Remark 4.6}  Rosenthal showed (unpublished) that any
weakly null sequence of indicator functions in a $C(K)$ space has a
subsequence which is an unconditional basic sequence by showing that
there is a weakly independent subsequence.\endremark

\remark{Remark 4.7}  Note that if $(x_n)$ is weakly independent and $
\Cal T$
is the associated tree then
$$\cup_k\left\{\left(n_1,n_2,\ldots
n_k\right):\left(\epsilon_1,\epsilon_
2,\ldots ,\epsilon_m\right)\in \Cal T\text{ where }\epsilon_{n_i}=
1\text{ for }i=1,2,\ldots ,k\right\}$$
$$=\Cal T\left(\left(\text{supp }x_n,\left(\text{supp }x_n\right)^{
-1}\right)\right),$$
the Boolean independence tree.
\vskip .15in
In order to write tree elements more efficiently let us introduce the
following notational conventions:
If $\overline x=\left(\epsilon_1,\epsilon_2,\ldots ,\epsilon_k\right
)$ and $\overline y=\left(\epsilon^{\prime}_1,\epsilon^{\prime}_2,
\ldots\epsilon^{\prime}_j\right)$ are two elements in
$\displaystyle{\cup^{\infty}_{n=1}S^n}$, for some set $S$, then
$\overline x+\overline y=\left(\epsilon_1,\epsilon_2,\ldots ,\epsilon_
k,\epsilon^{\prime}_1,\epsilon^{\prime}_2,\ldots ,\epsilon^{\prime}_
j\right)$, the concatenation of $\overline x$ and $\overline y$.  Let
$\overline e_0$ be the empty tuple, $\overline e_1=(-1)$,
and inductively define $\overline
e_{n+1}=\overline e_n+\overline e_1$,
$n=1,2,\ldots$ .  Also let $\overline e_{\omega}=(-1,-1,\ldots )$.

We will next introduce a derivation on a tree $\Cal T$ $\subset \Cal
D$.
Let
$$\delta^1(\Cal T)=\left\{\overline t\in \Cal T:
\text{ there are infinitely many }\overline
s\in\overline {\Cal T}\backslash \Cal T\text{ with }\overline
t<\overline
s\right\}$$
and if $\delta^{\alpha}(\Cal T)$ has been defined let $\delta^{\alpha
+1}(\Cal T)=\delta^1(\delta^{\alpha}(\Cal T))$.  If $\beta$ is a limit
ordinal let $\displaystyle{\delta^{\beta}(\Cal T)=\cap_{\alpha <\beta}
\delta^{\alpha}(\Cal T)}$.  Finally define
$\delta (\Cal T)=\inf\left\{\alpha :\delta^{\alpha}(\Cal
T)=\emptyset\right
\}$.\endremark

\proclaim{Proposition 4.8}  Suppose that $\Cal T$ is a weakly
independent
subtree of the dyadic tree with no infinite nodes, no infinite nodes in
$\overline {\Cal T}$ with infinitely many coordinates
equal to 1,  and every node of $
\Cal T$
is on some infinite branch.  Then $\delta (\Cal T$) determines the
homeomorphic
type of $\overline {\Cal T}$up to the number of points in the
last derived set, i.e.,
$\delta^{\alpha}(\Cal T)=\emptyset$ if and only if $\overline {\Cal
T}^{
(1+\alpha )}=\emptyset$.
\endproclaim

\demo{Proof}  First note that each element of $\overline {\Cal
T}\backslash
\Cal T$ is in the first
derived set of $\overline {\Cal T}$, and, in fact, $\overline {\Cal T}
\backslash \Cal T$ is the first derived set.  Now let
us set up a correspondence between the derived sets of $\overline {
\Cal T}$ and the
subtrees of $\Cal T$.  If $C$ is any closed subset of $\overline {
\Cal T}\backslash \Cal T$ then there is a tree
$\Cal T(C)$ with $\overline {\Cal T(C)}\backslash \Cal T(C)=C$.  Indeed
let $
\Cal T(C)$ be the set of all nodes $\overline x$ of
$\Cal T$ for which there is some element $c$ of $C$ below $\overline
x$.

As above let $C$ be a closed subset of $\overline {\Cal T}\backslash
\Cal T$.  We claim that
$\Cal T(C^{(1)})=\delta^1(\Cal T(C))$.  Suppose that $\overline x\in
\Cal T(C^{(1)})$.  Then $\overline x$ is above some
element $\overline c$ of $C^{(1)}$.  Say $\overline c=\overline
x+\overline
y+\overline e_{\omega}$ where $\overline y$ is possibly the empty
tuple.  Hence there are distinct points $\overline c_k$ in $C$
which converge to $\overline
c$
and therefore for large $k$, $\overline c_k=\overline y+\overline
z_k+\overline e_{\omega}$, where the $\overline z_k$ $^{\prime}s$ are
distinct.  This implies that $\overline x\in\delta^1(\Cal T(C))$.
Conversely suppose that
$\overline x\in\delta^1(\Cal T(C))$.  Then there is a
sequence of distinct points $
(\overline c_i)$ in
$\overline {\Cal T(C)}\backslash \Cal T(C)=C$ such that $\overline
c_i>\overline x$ for all $i$.  By passing to a subsequence
we may assume that $\overline c_i\longrightarrow\overline c\in C$.
Clearly $\overline
c\in C^{(1)}$ and $\overline c>\overline x$, therefore
$\overline x\in \Cal T$ $(C^{(1)})$.

Because $\overline {\Cal T}\backslash \Cal T=\overline {\Cal T}^{(
1)}$ and $\Cal T(\overline {\Cal T}\backslash \Cal T)=\Cal T$, we have
that
$$\delta^1(\Cal T)=\Cal T((\overline {\Cal T}\backslash \Cal T)^{(
1)})=\Cal T(\left[\overline {\Cal T}^{(1)}\right]^{(1)})=\Cal
T(\overline {
\Cal T}^{(2)}).$$
We claim that for every $\alpha <\omega_1$,
$$\delta^{\alpha}(\Cal T)=\Cal T(\left[\overline {\Cal T}^{(1)}\right
]^{(\alpha )})=\Cal T(\overline {\Cal T}^{(1+\alpha )}).$$
Indeed, if this is true for $\alpha$ we have by the previous claim that
$$\delta^{\alpha +1}(\Cal T)=\delta^1\left[\Cal T(\left[\overline {
\Cal T}^{(1)}\right]^{(\alpha )})\right]=\Cal T(\left[\left[\overline {
\Cal T}^{(1)}\right]^{(\alpha )}\right]^{(1)})=\Cal T(\left[\overline {
\Cal T}^{(1)}\right]^{(\alpha +1)}).$$
Also observe that if $(C_i)$ is a decreasing family of closed subsets
of
$\overline {\Cal T}\backslash \Cal T$ then\linebreak$\displaystyle{
\Cal T\left(\cap^{\infty}_{i=1}C_i\right)=\cap^{\infty}_{i=1}\Cal T
(C_i)}$.  Therefore if $\alpha_i\uparrow\alpha$ and
the claim is true for each $\alpha_i$ then
$$\delta^{\alpha}(\Cal T)=\cap\delta^{\alpha_i}(\Cal T)=\cap \Cal T
(\left[\overline {\Cal T}^{(1)}\right]^{(\alpha_i)})=\Cal T(\cap\left
[\left[\overline {\Cal T}^{(1)}\right]^{(\alpha_i)}\right])=\Cal
T(\left
[\overline {\Cal T}^{(1)}\right]^{(\alpha )}),$$
establishing the claim.  Hence $\delta^{\alpha}(\Cal T)=\emptyset$
if and only if $\overline {
\Cal T}^{(1+\alpha )}=\emptyset$.\qed\enddemo

Let us now return to the Schreier sequence and examine the
associated tree $\Cal S_1$.  Let
$$\align \Cal S_0=\cup^{\infty}_{n=0}\big\{\left(\epsilon_1,\epsilon_
2,\ldots ,\epsilon_n\right):&\epsilon_i=1\text{ for at most one }i
,\\
&\epsilon_i=-1\text{ otherwise, }1\le i\le n\big\}.\endalign$$
To define $\Cal S_1$ we need introduce the extension of one tree by
another
tree.  If $\Cal T$ and $\Cal S$ are trees on the same set let $\Cal
T\boxplus
\Cal S$ denote
$\left\{\overline x+\overline y:\overline x\in \Cal T\text{ and
}\overline
e_n+\overline y\in \Cal S\text{ where $n$ is the length}\right.$
\linebreak $\displaystyle{\text{of }
\left.\overline
x\right\}\cup \Cal T}$.  If $(\Cal T_i)$ is
a sequence of trees, we inductively define
$\displaystyle{\boxplus^n\Cal T_i=\left[\boxplus^{n-1}\Cal T_i\right
]\boxplus \Cal T_n}$.  If $n=0$, let
$\displaystyle{\boxplus^n\Cal T_i=\left\{\overline
e_j:j=1,2,\ldots\right
\}}$.  Also we will use the notation
$\Cal L(\Cal T,n)$ for the subtree which is the union of
$\left\{\overline
e_j:j=0,1,2,\ldots ,n-1\right\}$ and
the set of nodes of $\Cal T$ equal to or below $\overline e_{n-1}+
(1)$.  In particular
$\Cal L(\Cal T,n)\subset\left\{\overline e_j:j=1,2,\ldots ,n-1\right
\}$$\boxplus \Cal T$.  We need a way of forming a new
tree $\Cal T$ out of an infinite sequence of trees $\Cal T_i$, $i=
1,2,\ldots$ , on $\{-1,1\}$.

Define $\sum^{\infty}_{i=1}\Cal T_i=\cup^{\infty}_{i
=1}\Cal L(\Cal T_i,i)$
Clearly the resulting tree will depend on the order of the $\Cal T_
i$ $^{\prime}s$.

We claim that the tree associated to the Schreier sequence is
$$\Cal S_1=\sum^{\infty}_{i=1}\boxplus^i\Cal S_0.$$
Indeed, the Schreier sequence $(x_n)$ is defined by the property that
if
$k,N\in \Bbb N$ and $k\le n_1<n_2<\ldots <n_k\le N$, then
$$\bigcap^k_{i=1}\text{supp }x_{n_i}\cap
\bigcap\Sb n\neq n_i\\n<N\endSb \text{(supp
}x_n)^{-1}\neq\emptyset,$$
and these are the only non-empty intersections.  Note that this is
equivalent to saying that if $m\le k$ and $k=n_1<n_2<\ldots <n_m$,
$\displaystyle{\cap^m_{i=1}\text{supp }x_{n_i}\cap\cap\Sb n\neq n_i\\
n<N \endSb(\text{supp }x_n)^{-1}\neq\emptyset}$.
Hence the tree associated to $(x_n)$ contains for each $k\in \Bbb N$,
exactly those nodes of the dyadic tree of the form
$\overline e_{k-1}+\overline x$, where $\overline x$ has at most $
k$ coordinates equal to 1.
It is easy to see that $\displaystyle{\boxplus^k\Cal S_0}$ is
exactly the nodes with
at most $k$ coordinates equal to 1.  Thus $\Cal S_1$ is the required
tree.

To see what the topological type of $\overline {\Cal S}_1$ is
we need only compute the
order of $\displaystyle{\delta\left[\boxplus^n\Cal S_0\right]}.$
Clearly $
\delta (\Cal S_0)=2$.  A
straight-forward induction argument shows that
$\displaystyle{\delta^n(\left[\boxplus^n\Cal S_0\right])=}$
\linebreak$\displaystyle{
\delta^n(\Cal L(\boxplus^n\Cal S_0,n)=\{e_{\omega}\}}$.  It then
follows easily
that $\delta (\Cal S_1)=\omega +1$.  Therefore $\overline {\Cal S}_
1$ has exactly $\omega$ non-empty derived sets.
 Actually it is not hard to see that it is homeomorphic to $\omega^{
\omega}$.

Now we are ready to generalize the Schreier example.  Our goal is to
build for each $\alpha <\omega_1$ a weakly null sequence of indicator
functions
on a compact metric space (homeomorphic to $\omega^{\omega^{\alpha}}$)
with oscillation
index at least $\omega^{\alpha}$.  The sequence will
also be a 1-unconditional basic
sequence.

We will begin by defining the sequences in terms of subsets of $\Bbb N$
as
we did with the Schreier sequence.  $\Cal F_1$ has been defined.
Suppose
that $\Cal F_{\beta}$ has been defined for all $\beta <\alpha$.  Let $
\alpha_i$ $=\alpha -1$ if $\alpha$ is not a
limit ordinal and $\alpha_i\uparrow\alpha$ if $\alpha$
is a limit ordinal.  Define
$$\Cal F_{\alpha}=\cup^{\infty}_{n=1}\left\{\cup^n_{i=1}F_i:F_i\in
\Cal F_{\alpha_i},\text{ for }i\le n,\text{ }k\le F_1<\ldots <F_i<
F_{i+1}<\ldots F_n\right\}$$
where the notation $k\le F_1<\ldots F_i<F_{i+1}<\ldots F_n$ means that
if $
F_i$ and
$F_j$ are nonempty and $i<j$, then $k\le\min F_i$ and $\max F_i<\min
F_j$.

There is some ambiguity here but it will not be of any significance
as long as the choice of the sequence $(\alpha_i)$ is fixed for each
limit
ordinal.  For each $\alpha <\omega_1$ let $(x^{\alpha}_n$)
denote the standard sequence of
indicator functions on $\Cal F_{\alpha},$ i.e., $x^{\alpha}_n(F)=1$ if
$
n\in F$ and 0 otherwise.

\proclaim{Proposition 4.9}  For each $\alpha <\omega_1$, $\Cal F_{
\alpha}$ is a countable
compact metric space under the topology induced by identifying $\Cal
F_{
\alpha}$
with $\left\{1_F:F\in \Cal F_{\alpha}\right\}$ under the
topology of pointwise convergence and
$\left(x^{\alpha}_n\right)$ is a weakly null 1-suppression
unconditional basic sequence in
$C(\Cal F_{\alpha})$.\endproclaim

\demo{Proof}  It is sufficient to verify that $\Cal F_{\alpha}$
satisfies the
hypothesis of Proposition 4.2.    Property ii) is obviously inherited
by each $\Cal F_{\alpha}$ from $\Cal F_1$.  For i)
use induction and note that if
$G\subset F_1\cup F_2\cup\ldots\cup F_n\in \Cal F_{\alpha}$,
as in the definition, then $
G\cap F_i\in \Cal F_{\alpha_i}$, for
each $i$.  Hence  $G=(G\cap F_1)\cup (G\cap F_2)\cup\ldots\cup (G\cap
F_n)$ $\in \Cal F_{\alpha}$.  Finally for
iii) suppose that iii) holds for all $\beta <\alpha$, and that for each
k,
$\displaystyle{F_k=\cup^{n_k}_{j=1}F_{kj}}$ is an element of $\Cal F_{
\alpha}$, as in the definition.
If $F_k$ converges to a non-empty set $F$, let $n=\min F$.  Then $
n\in F_k$ for
all large $k$ and thus $n_k\le n$.  We may assume by passing to
subsequences that for each $j\le n$, $F_{kj}$ converges to some $F_
j$.  Hence,
by induction, because we have that $F_{kj}\in \Cal F_{\alpha j}$, for
all $
k$, $F_j\in \Cal F_{\alpha j}$.
the other properties are obvious and thus $F=\cup F_j\in \Cal
F_{\alpha}$, as
claimed.\qed\enddemo

Now we want to consider the size of the underlying topological
space.  First we will compute the size directly using the families
$\Cal F_{\alpha}$.  For each $\alpha <\omega_1$ and $k$, $n\in \Bbb N$
with $
n\le k$ let
$$\Cal F_{\alpha ,n,k}=\left\{\cup^n_{i=1}F_i:F_i\in \Cal F_{\alpha_
i},k\le F_1<\ldots <F_i<F_{i+1}<\ldots <F_n,\text{ and }F_n\in \Cal F_{
\alpha_n}\right\},$$
where $\alpha_i=\alpha -1$ if $\alpha$ is not a limit ordinal and $
\alpha_i\uparrow\alpha$ if $\alpha$ is a limit
ordinal.

\proclaim{Proposition 4.10}  For each $\alpha <\omega_1$, $\Cal F^{
(\omega^{\alpha})}_{\alpha}=\{\emptyset \}$.\endproclaim

\demo{Proof}  The result will follow from

CLAIM:  If $\rho\le\omega^{\alpha_n}$, then
$$\Cal F^{(\rho )}_{\alpha ,n,k}=\left\{\cup^n_{i=1}F_i:F_i\in \Cal F_{
\alpha_i},k\le F_1<\ldots <F_i<F_{i+1}<\ldots <F_n,\text{ and }F_n
\in \Cal F_{\alpha_n}^{(\rho )}\right\}$$
The proof of the claim is by induction on $\alpha ,n,$ and $\rho$.

If $\alpha =0,$ $\Cal F_0=\left\{\left\{n\right\}:n\in \Bbb N\right
\}\cup \{\emptyset \}$.  Clearly $\Cal F^{(1)}_0=\{\emptyset \}$.

Let $\alpha =1$ and $k\ge n\ge 1$.  Then if $G\in \Cal F^{(1)}_{1,
n,k}$, there is a non-trivial
sequence $(G_m)$ in $\Cal F_{1,n,k}$ which converges to $F$.
Observe that we need
only show that card $G<n$.  However this is obvious because card
$G_m\le n$ for all $m$ and some portion of the $G_m$ $^{\prime}s$ must
go to $
\infty$.
Conversely, if $G\in \Cal F_{1,n,k}$ and card $G<n$, then $G\cup \{
m\}\in \Cal F_{1,n}$ for $m=k,$
$k+1,\ldots.$  Hence $G\in \Cal F^{(1)}_{1,n,k}$.  Finally note that
each $
\Cal F_{1,n,k}$ is closed,
$\Cal F^{(1)}_{1,n,k}=\Cal F_{1,n-1,k}$, and if $G_m\in \Cal F_{1,
m,m}$ then $G_m\longrightarrow\emptyset$.  Therefore $\Cal F^{(\omega
)}_1=\{\emptyset \}$.

Now assume the result holds for all $\gamma <\alpha$.  Fix $\rho <
\omega^{\alpha_n}$ and $k\ge n\ge 1$
and assume that the result has been proved for $\rho$.  (Note that it
always holds for $\rho =0.)$  Let $\Cal F=\Cal F^{(\rho )}_{\alpha
,n,k}$.  Suppose that $(G_m)$ is a
non-trivial sequence in $\Cal F$ which converges to $G$.  Suppose that
$\displaystyle{G_m=\cup^n_{i=1}G_{m,i}}$.  By passing to a subsequence
if
necessary we may assume that for each $i$, $G_{m,i}\longrightarrow
G_i$.  If $(G_{m,n})$ is a
non-trivial sequence, then $G_n\in \Cal F^{(\rho +1)}_{\alpha_n}$ and
thus
$\displaystyle{G=\cup^n_{i=1}G_i}$ where $G_i\in \Cal F_{\alpha_i}$ for
$
i=1,2,\ldots ,n$ and
$G_n\in \Cal F^{(\rho +1)}_{\alpha_n}$, that is,
$$G\in\left\{\cup^n_{i=1}F_i:F_i\in \Cal F_{\alpha_i},k\le F_1<\ldots
<F_i<F_{i+1}<\ldots <F_n,\text{ and }F_n\in \Cal F^{(\rho +1)}_{\alpha_
n}\right\},$$
If $(G_{m,n})$ is trivial (eventually constant) then $(G_{i,n})$ is
eventually
constant for all $i$.  However this contradicts the non-triviality of
the
sequence $(G_m)$.

Conversely if
$$G\in\left\{\cup^n_{i=1}F_i:F_i\in \Cal F_{\alpha_i},k\le F_1<\ldots
<F_i<F_{i+1}<\ldots <F_n,\text{ and }F_n\in \Cal F^{(\rho +1)}_{\alpha_
n}\right\},$$
Then there is a non-trivial sequence $(G_m)$ in $\Cal F^{(\rho )}_{
\alpha_n}$ which converges to
$F_n$.  Clearly we may assume that $\displaystyle{\min G_m>\max \cup^{
n-1}_{i=1}F_i}$
for all $m$.  Then $\displaystyle{G^{\prime}_m=\cup^{n-1}_{i=1}F_i
\cup G_m\in \Cal F^{(\rho )}_{\alpha ,n,k}}$ for all $m$
and $(G^{\prime}_m)$ converges to $G$.  Therefore $G\in \Cal F^{(\rho
+1)}_{\alpha ,n,k}$.

Clearly if $\rho_j\uparrow\rho$, $\Cal F^{(\rho )}_{\alpha ,n,k}=$
$$\cap^{\infty}_{j=1}\left\{\cup^n_{i=1}F_i:F_i\in \Cal F_{\alpha_
i},k\le F_1<\ldots <F_i<F_{i+1}<\ldots <F_n,\text{ and }F_n\in \Cal F^{
(\rho_j)}_{\alpha_n}\right\}$$
$$=\left\{\cup^n_{i=1}F_i:F_i\in \Cal F_{\alpha_i},k\le F_1<\ldots
<F_i<F_{i+1}<\ldots <F_n,\text{ and }F_n\in \Cal F^{(\rho )}_{\alpha_
n}\right\}.$$
Thus the formula holds for all ordinals $\rho\le\omega^{\alpha_n}$.

Finally observe that by induction we have that $\Cal
F^{(\omega^{\alpha_
n})}_{\alpha ,n,k}=\Cal F_{\alpha ,n-1,k}$ and
hence that\linebreak\ $\Cal F^{(\omega^{\alpha_n}+\omega^{\alpha_{
n-1}}+\ldots +\omega^{\alpha_1})}_{\alpha ,n,n}=\{\emptyset \}$.  To
see that $
\Cal F^{(\omega^{\alpha})}_{\alpha}=\{\emptyset \}$ note
that $\displaystyle{\Cal F_{\alpha}=\cup^{\infty}_{n=1}\Cal F_{\alpha
,n,n}}$ and that if $G_n\in \Cal F_{\alpha ,n,n}$ for all $n$
then $G_n\longrightarrow\emptyset$.\qed\enddemo

Next we will prove the result again but using trees.  First we will
translate the construction into the tree representation.

Observe that there is a simple correspondence between
$$\left\{\cup^n_{i=1}F_i:F_i\in \Cal F_{\alpha_i},k\le F_1<\ldots
<F_i<F_{i+1}<\ldots <F_n\right\}$$
and $\displaystyle{\boxplus_{i=1}^n\Cal S_{\alpha_i}}$, where $\Cal S_{
\alpha_i}$ is the tree corresponding to
$\Cal F_{\alpha_i}$.  Indeed,
$$\bigcap_{n\in F_i}\text{supp }
x^{\alpha}_n\cap\bigcap\Sb n\notin F_i\\
n\le m \endSb(\text{supp }x_n^{\alpha})^{-1}\neq\emptyset,$$
 for all $m$, if and only if $F_i\in \Cal F_{\alpha}$.  Thus if
 $\overline
y_i=\left(\epsilon_{m_i+1},\epsilon_{m_i+2},\ldots ,\epsilon_{m_{i
+1}}\right)$,
where $m_i=\min F_i-1$, and $\epsilon_j=1$ if $j\in F_i$, $\epsilon_
j=-1$, otherwise, then we
have $\overline y_1+\overline y_2+\ldots +\overline y_n+\overline
e_m$ in $\displaystyle{\boxplus^n\Cal S_{\alpha}}$, for all $m$.
Conversely if we have $\overline y_1+\overline y_2+\ldots +\overline
y_n$ in $\displaystyle{\boxplus^n\Cal S_{\alpha}}$, let $m_i$
be the sum of the lengths of the $\overline y_j$ $^{\prime}s$, $j=
1,2,\ldots ,i-1$, and
$F_i=\left\{m_i+k:\epsilon_k=1\right\}$ for $\overline
y_i=\left(\epsilon_
1,\epsilon_2,\ldots ,\epsilon_{m_{i+1}-m_i}\right)$.  Then $F_i\in
\Cal F_{\alpha},$ and
$\max F_i<\min F_{i+1},$ for $i=1,2,\ldots ,n$.

Next we need to take care of the $n\le\min F_1$ condition.  Observe
that
we could define
$$\Cal F_{\alpha}=\cup^{\infty}_{n=1}\left\{\cup^{\infty}_{n=1}F_i
:F_i\in \Cal F_{\alpha_i},n\le F_1<\ldots <F_i<F_{i+1}<\ldots
<F_n\text{ and }
n\in F_1\right\}$$
and the set would remain the same.  We claim that the tree
associated to $x^{\alpha +1}_n$ is
$$\Cal S_{\alpha +1}=\sum^{\infty}_{i=1}\boxplus^i\Cal S_{\alpha}.$$
We have already shown above that $\displaystyle{\boxplus^i\Cal S_{
\alpha}}$ corresponds
to the union of $i$ ordered sets from $\Cal F_{\alpha}$.
Now note that the nodes
corresponding to $\displaystyle{\boxplus^i\Cal S_{\alpha}}$ in
the sum are all below or
equal to $\overline e_{i-1}+(1)$, and thus correspond
exactly to the unions of $
i$
orders sets from $\Cal F_{\alpha}$ with smallest
element of the first set equal to
$i$.  Thus $\Cal S_{\alpha +1}$ is the correct tree.
A similar argument shows that
for a limit ordinal $\alpha$
$$\Cal S_{\alpha}=\sum^{\infty}_{i=1}\boxplus_{j=1}^i\Cal S_{\alpha_
j}$$
where $\alpha_i$ is the defining sequence for $\Cal F_{\alpha}$.

Now that we have the trees $\Cal S_{\alpha}$, $\alpha <\omega_1$,
we can determine the
underlying topological spaces by using Proposition 4.7.

\demo{Proof of Proposition 4.10}  Inductively assume that
$\delta (\Cal L(\Cal S_{\beta},j))$ is $\omega^{\beta}$, for all $
\beta <\alpha$ and $j\in \Bbb N$.  Let $\alpha_i=\alpha -1$ if $\alpha$
is a
successor and $\alpha_i\uparrow\alpha$ otherwise.\enddemo
$$\text{CLAIM:  }\delta (\Cal L(\boxplus_{i=1}^n\Cal S_{\alpha_i},
j))=\omega^{\alpha_n}+\omega^{\alpha_{n-1}}+\ldots +\omega^{\alpha_
1}+1\text{ for all }j\in \Bbb N.$$
Indeed if $\overline x\in \Cal L\left(\boxplus_{i=1}^{n-1}\Cal S_{
\alpha_i},j\right)$, $\overline x=$ $\overline y_1+\overline y_2+\ldots
+\overline y_{n-1}$ and the length of $\overline x$
is $k\ge j$, then
$$\left\{\overline e_k+\overline y:\overline x+\overline y\in \Cal
L\left
(\boxplus_{i=1}^n\Cal
S_{\alpha_i},k+1\right)\right\}\cup\left\{\overline
e_i:i=1,2,\ldots ,k\right\}
\supset \Cal L\left(\Cal S_{\alpha_n},k+1\right)$$
Thus
$$\overline x\in\delta^{\omega^{\alpha_n}}\left(\Cal L\left(\boxplus^
n\Cal S_{\alpha_i},j\right)\right),$$
by the inductive hypothesis.  Because $k$ is arbitrary the claim
follows by induction on $n$.  (Obviously
$\displaystyle{\delta\left(\Cal L\left(\boxplus_{i=1}^n\Cal S_{\alpha_
i},j\right)\right)}$ $\le\omega^{\alpha_n}+\omega^{\alpha_{n-1}}+\ldots
+\omega^{\alpha_1}+1$.)

It now is easy to see that $\delta\left(\Cal L\left(\Cal S_{\alpha}
,j\right)\right)=\omega^{\alpha}+1$, because the order is
larger than $\lambda_n=\omega^{\alpha_n}+\omega^{\alpha_{n-1}}+\ldots
+\omega^{\alpha_1},$ for all $n$ and the elements $e_i$,
$i\in \Bbb N$, are in $\delta^{\lambda_n}\big(\Cal L\left(\boxplus^
n\Cal S_{\alpha_i},j\right)$ for all $n$.\qed

Our next task is to compute the oscillation index of $\left(x^{\alpha}_
n\right)$.  As before
we will compute this in two ways first by using the family $\Cal F_{
\alpha}$ and
then by using trees.

\proclaim{Proposition 4.11}  $\Cal O^{\omega^{\alpha}}\left(\epsilon
,\left(x^{\alpha}_n\right),\Cal F_{\alpha}\right)\neq\emptyset$, for
all $
\alpha <\omega_1,$ $\epsilon <1$.
\endproclaim

\demo{Proof}  We will show that $\Cal O^{\lambda}\left(\epsilon ,\left
(x^{\alpha}_n\right),\Cal F_{\alpha}\right)
=\Cal F^{(\lambda )}_{\alpha}$ for all $
\lambda$.  In
view of Proposition 4.10, it is sufficient to show that if $F\in \Cal
F^{
(\lambda +1)}_{\alpha}$
then there is a $N\in \Bbb N$ such that for all $n\ge N$, $F\cup \{
n\}\in \Cal F^{(\lambda )}_{\alpha}$.  We will
use induction on $\alpha$ and $\lambda$.

If $\alpha =1$, then for any $\lambda\in \Bbb N$ the claim
is immediate from the
definition of $\Cal F_1$ and Proposition 4.10.

Now assume that $\alpha >1$ and that the claim is true for all $\beta
<\alpha$ and
all $\lambda$.  If $F\in \Cal F^{(1)}_{\alpha}$ then
either $F=\emptyset$ and the claim is obvious or
$F\in \Cal F^{(1)}_{\alpha ,k,k}$ for some $k\in \Bbb N$.  By
Proposition 4.10,
$$\Cal F^{(1)}_{\alpha ,k,k}=\left\{\cup^k_{i=1}F_i:F_i\in \Cal F_{
\alpha_i},k\le F_1<\ldots <F_i<F_{i+1}<\ldots <F_k,\text{ and }F_k
\in \Cal F^{(1)}_{\alpha_k}\right\}.$$
Suppose that $\displaystyle{F=\cup^k_{i=1}F_i}$ as above.  By the
inductive
hypothesis $F_k\cup \{n\}\in \Cal F_{\alpha_k}$ for all $n\ge N$, for
some $
N\in \Bbb N$, and thus for
$n>\max \{N\}\cup F_k$, $\displaystyle{\cup^{k-1}_{i=1}F_i\cup F_k
\cup \{n\}\in \Cal F^{(0)}_{\alpha ,k,k}}$.

Next assume the claim for all $\gamma\le\lambda$ and let $F\in \Cal F^{
(\lambda +1)}_{\alpha}$.  If $F=\emptyset$ then
$\{n\}\in \Cal F^{(\lambda )}_{\alpha}$ for all sufficiently large $
n$.  Otherwise $F\in \Cal F^{(\lambda +1)}_{\alpha ,k,k}$ for some
$k$.  Let $\lambda =\omega^{\alpha_k}+\omega^{\alpha_{k-1}}+\ldots
+\omega^{\alpha_j}+\rho$ for some $j\le k+1$ and $\rho
<\omega^{\alpha_{
j-1}}$.
By the claim in the proof of Proposition 4.10
$$\multline\Cal F^{(\lambda +1)}_{\alpha ,k,k}=\Cal F^{(\rho +1)}_{\alpha ,
j-1,k}\\
=\left\{\cup^{j-1}_{i=1}F_i:F_i\in \Cal F_{\alpha_i},k\le F_1<\ldots
<F_i<F_{i+1}<\ldots <F_{j-1},\text{ and }F_{j-1}\in \Cal F^{(\rho
+1)}_{\alpha_{j-1}}\right\}.\endmultline$$
Suppose that $\displaystyle{F=\cup^{j-1}_{i=1}F_i}$ as above.  By the
inductive
hypothesis $F_{j-1}\cup \{n\}\in \Cal F_{\alpha_{j-1}}$ for all $n
\ge N$, for some $N\in \Bbb N$, and thus
for $n>\max \{N\}\cup F_{j-1},$ \linebreak
$\displaystyle{\cup^{j-2}_{i=1}F_i
\cup F_{j-1}\cup \{n\}\in \Cal F^{(\rho )}_{\alpha
,j-1,k}}$.\qed\enddemo

In the proof above we established that we can always use the special
sequence\linebreak\ $\left(F\cup \{n\}\right)_{n\ge N}$ to reach a set
$
F$ from a smaller derived set.
The next definition describes this same property for the associated
tree.

\definition{Definition}  Let $\Cal T$ be a tree on $\{-1,1\}$ with
no nodes in $\overline {
\Cal T}$
with infinitely many coordinates equal to 1.  We will say that $\Cal T$
has
{\it property\/} FB (fully branching) if $\overline x+$ $\overline
e_{\omega}\in\overline {\Cal T}$ implies that there is an
$N\in \Bbb N$ such that either

\itemitem{i)}  $\overline x+\overline e_j+(1)\in \Cal T$ for all $
j\ge N$
or
\itemitem{ii)}  $\overline x+\overline e_j+(1)\notin \Cal T$ for all $
j\ge N$.
\enddefinition

\proclaim{Lemma 4.12}  Suppose that $\Cal T$ and $\Cal U$
are weakly independent
trees with property FB and all nodes of $\Cal T$ and $\Cal U$ are on
branches
with limit of the form $x+e_{\omega}$, then $\Cal T\boxplus \Cal U$
has property FB, and for all
$n$, $\Cal L\left(\Cal T,n\right)$ has property FB.
\endproclaim

\demo{Proof}  Suppose that $\overline x+\overline
e_{\omega}\in\overline {
\Cal T}$.  Because $\Cal T$ $\boxplus \Cal U\supset \Cal T$, if i)
occurs
in $\Cal T$, the same is true in $\Cal T\boxplus \Cal U$.  If ii)
occurs in $
\Cal T$ but not in $\Cal T\boxplus \Cal U$, then
there is a sequence of incomparable nodes of the form $\overline e_
n+\overline y_j$ in $\Cal U$
where $n$ is greater than the length of $\overline x$ and does not
depend on $
j$ and
$\overline y_j$ has at least one coordinate equal to one.  By passing
to a
subsequence we may assume that $\overline e_n+\overline y_j$
converges to $\overline
e_n+\overline z+\overline e_{\omega}$.
The assumption that the $\overline e_n+\overline y_j$ $^{\prime}s$
are incomparable guarantees that
$(\overline e_n+\overline z)+\overline e_{\omega}$ does not
satisfy ii).  Because $
\Cal U$ has property FB there is
an $N$ such that $\overline e_n+\overline z+\overline e_j+(1)\in \Cal
U$
for all $
j>N$.  Because $\Cal U$ is weakly
independent this implies that $\overline e_k+\overline e_j+(1)\in
\Cal U$ for all $j>N$, where $k$
equals $n$ plus the length of $\overline z$.  Hence $\overline x$ $
+\overline e_m+(1)\in \Cal T\boxplus \Cal U$ for all
$m>N+k,$ i.e., $\overline x+\overline e_{\omega}$ satisfies i).

If $\overline x+\overline y+\overline e_{\omega}\in\overline {\Cal
T\boxplus
\Cal U}\backslash\overline {\Cal T}$, where $\overline e_n+\overline
y\in \Cal U$, $n$ is the length of $\overline x$ and $\overline x\in
\Cal T$,
then $\overline x+\overline y+\overline e_{\omega}$
will satisfy i), respectively ii), if $\overline
e_n+\overline y+\overline e_{\omega}$ satisfies
i), respectively ii).

The second assertion is obvious.\qed\enddemo

\proclaim{Proposition 4.13}  Suppose that $(x_n)$ is a weakly
independent sequence of continuous indicator functions on a compact
metric space $K$ which converge pointwise to 0 and that the
associated tree $\Cal T$ has property FB.  Then for any $\epsilon
<1$, and $\alpha <\omega_1$,
$\Cal O^{\alpha}\left(\epsilon ,\left(x_n\right),K\right)\neq\emptyset$
if and only if $\overline {
\Cal T}^{(1+\alpha )}\neq\emptyset$.\endproclaim

\demo{Proof}  For each $n$ and $(\epsilon_i)$ $\in\overline {\Cal T}$ $
\backslash \Cal T$, let $\hat {x}_n\left(\left(\epsilon_i\right)\right
)=1$ if $\epsilon_n=1$, and
0, otherwise.  In this way we have defined a sequence of indicator
functions on $Q=\overline {\Cal T}\backslash \Cal T$
with span isometric to the span of $
(x_n)$.
(Actually $\overline {\Cal T}\backslash \Cal T$ is
homeomorphic to the natural quotient of $ K$
determined by the $x_n$ $^{\prime}s$.)

Clearly $\Cal O^{\alpha}\left(\epsilon ,\left(x_n\right),K\right)\neq
\emptyset$ if and only if $\Cal O^{\alpha}\left(\epsilon
,\left(\hat x_n\right),Q\right)\neq\emptyset$.  Now observe that
$(\epsilon_i)+\overline {\epsilon}_{\omega}\in \Cal O^ 1\left(\epsilon
,\left(\hat x_n\right),Q\right)$ if and only if there is an $ N\in
\Bbb N$ such that for all $j\ge N$, $(\epsilon_j)+\overline e_j+(1)\in
\Cal T$.  (Use the weak independence of $\left(\hat x_n\right).)$
Clearly this latter condition implies that $(\epsilon_i)+\overline
e_{\omega}\in\overline {\Cal T}^{(2)}=Q^{( 1)}$.  Conversely if
$(\epsilon_1)\in Q^{(1)}$, then weak independence and property FB imply
that $(\epsilon_i)+\overline e_j+(1)\in \Cal T$, for all $ j\ge N$, for
some $N\in \Bbb N$.  Finally note that weak independence and property FB
are inherited by $\Cal T(Q^{(1)})$ which is the tree associated to $\left
(\hat x_{n|_{Q^{^{(1)}}}}\right)$.  Also if $\alpha_k\uparrow\alpha ,$
$\cap \Cal T\left(Q^{(\alpha_n)}\right) =\Cal T\left(Q^{(\alpha )}\right)$
and it is straight-forward to check that weak independence and property
FB are inherited by the intersection.  Thus transfinite induction may
be used to complete the proof.\qed\enddemo

Proposition 4.11 follows as a corollary of this result, that is,
$\Cal O^{\omega^{\alpha}}\left(\epsilon ,\left(x^{\alpha}_n\right)
,K\right)\neq\emptyset$, for all $\alpha <\omega_1$, $\epsilon <1$.

Because the underlying space for $\left(x^{\alpha}_n\right)$ is
$\omega^{
\omega^{\alpha}}$, $\omega^{\alpha}$ is the maximal
possible oscillation index.  Also note that because each $\left(x^{
\alpha}_n\right)$ is an
unconditional basic sequence, all of the spaces are isomorphic to
complemented subspaces of the Pelczynski universal space $U_1$,
[L-T,I,p92].

\head 5.  Reflexive spaces with large oscillation index\endhead

In the previous section we constructed weakly null sequences with
oscillation index $\omega^{\alpha}.$  Because these sequences were in $
C\left(\Cal F_{\alpha}\right)$ it
follows that the span of $\left(x^{\alpha}_n\right)$ contains $c_0$ and
thus is not reflexive.
In this section we will explore an idea of E. Odell for constructing
Tsirelson-like spaces with large oscillation index.

To define these spaces we begin with the space
$\left[x^{\alpha}_n\right
]$ in place of $c_0$
in the Tsirelson construction, [C-S,p 14].  Suppose that $x=\sum a_
nt^{\alpha}_n$
where $\left(t^{\alpha}_n\right)$ is the unit vector basis of the
space of sequences with
only finitely many nonzero coordinates.  Let
$$||x||_0=||\sum a_nx^{\alpha}_n||$$
and inductively define
$$||x||_{m+1}=\max\left\{||x||_m,2^{-1}\max_{\{p_i\}\in \Cal F_1}\sum^
k_{i=1}||\sum^{p_{i+1}-1}_{n=p_i}a_nt^{\alpha}_n||_m\right\},$$
where $k$ is the cardinality of $\{p_i\}$.  Let $||x||=\lim ||x||_
m$.  Let $T_{\alpha}$ be the
completion of span $\{t^{\alpha}_n$ $\}$ under $||\cdot ||$.
Observe that Tsirelson space is
$T_0$ in this construction, i.e., let $\Cal F_0=\left\{\left\{n\right
\}:n\in \Bbb N\right\}$, then $\left(x^0_n\right)=\left(1_{\left\{\left
\{n\right\}\right\}}\right)$
is equivalent to the $c_0$ basis.  Also observe that for each $\alpha$
the norm
on $T_{\alpha}$ satisfies
$$||x||=\max\left\{||x||_0,2^{-1}\max_{\{p_i\}\in \Cal F_1}\sum^k_{
i=1}||\sum^{p_{i+1}-1}_{n=p_i}a_nt^{\alpha}_n||\right\}.\tag *$$

\proclaim{Proposition 5.1}  For each $\alpha <\omega_1$,
$T_{\alpha}$ is a reflexive Banach
space with unconditional basis with no subspace isomorphic to $\ell^
p$,
$1\le p<\infty$, or $c_0$.  Moreover $\Cal O^{\omega^{\alpha}}\left
(\epsilon
,\left(t^{\alpha}_n\right),B_{\left(T_{\alpha}\right)^{*}}\right
)\neq\emptyset$ for all $\epsilon <1$.
\endproclaim

\demo{Proof}  Fix $\alpha <\omega_1$.  Clearly $||x||_0\ge ||x||_{
c_0}$.  Therefore
$||x||\ge ||x||_T$ for all $x\in T_{\alpha}$ where $||\cdot ||_T$
denotes the norm on Tsirelson
space.  (This is $T_0$ above.)  It follows from $(^{*})$ that if $
\{u_i\}$ is a
sequence of $k$ normalized blocks of the basis in $T_{\alpha}$ with
support
beyond $k$ that
$$\sum^k_{i=1}|c_i|\ge ||\sum^k_{i=1}c_iu_i||_T\ge 2^{-1}\sum^k_{i
=1}|c_i|,$$
for all choices of scalars $(c_i)$.  Therefore $T_{\alpha}$ does not
contain $
c_0$ or
$\ell^p$, for any $p>1$.

To see that $\ell^1$ is not isomorphic to a
subspace of $T_{\alpha}$ we need only
examine the proof that $T$ does not contain $\ell^1$ as given in [C-S,p
17].
The only properties of $T$ that are used in the argument are that the
norm satisfies equation $(^{*})$ and that the $c_0$ norm of a long
average is
small.  The proof then for $T$ will carry over to $T_{\alpha}$ provided
we use
a sequence with small $||\cdot ||_0$.  To do this let $(u_i)$ be a
normalized
block basis of $(t^{\alpha}_n)$ such that
$$\sum |a_i|\ge ||\sum a_iu_i||\ge\frac 89\sum |a_i|.$$
Because $\left(x^{\alpha}_n\right)$ is weakly null in $C(\Cal
F_{\alpha}
)$, so is $(u_i)$ and thus there is a
sequence of disjoint convex combinations of $(u_i)$,
$$y_j=\sum_{i\in E_j}\lambda_iu_i,$$
where $\displaystyle{E_1<E_2<\ldots ,\sum_{i\in E_j}\lambda_i=1}$, and
$
\lambda_i\ge 0$ for all $i$,
such that $||y_j||_0\le 2^{-j}$ for $j=0,1,\ldots$ .  It follows then
that
$$\sum |a_i|\ge ||\sum a_iy_i||\ge\frac 89\sum |a_i|,$$
and
$$||y_0+r^{-1}(y_1+y_2+\ldots +y_r)||_0\le 1+r^{-1}.$$
Using these $y_i$ $^{\prime}s$ the remainder of the proof carries over
without
change to $T_{\alpha}$.

Obviously $T_{\alpha}$ has an unconditional basis and thus by a result
of
James [L-T,I,p.97], $T_{\alpha}$ is reflexive.

Finally to see that the oscillation index is large we use the
observation that the operator $S$ from $T_{\alpha}$ to
$\left[x^{\alpha}_
n\right]$ defined by
$S(\sum c_nt_n)=\sum c_nx^{\alpha}_n$ is bounded by 1 and thus
$\Cal O^{\lambda}\left(\epsilon ,\left(t^{\alpha}_n\right),B_{T^{*}_{
\alpha}}\right)\supset S^{*}\Cal O^{\lambda}\left(\epsilon ,\left(
x^{\alpha}_n\right),\Cal F_{\alpha}\right)\neq\emptyset$, for every $
\lambda\le\omega^{\alpha}$, by Lemma
2.2.\qed\enddemo

While these space $T_{\alpha}$ share important properties with $T$ let
us note
that they do not possess the property that every block basis
dominates a subsequence of the basis.  In particular

\proclaim{Proposition 5.2}  For every $\alpha <\omega_1$ there is a
block basis
$(u_i)$ of $\left(t^{\alpha}_n\right)$ and an increasing sequence of
integers $
(k_i)$ such that
$$||\sum a_it_i||_T\le ||\sum a_iu_i||\le 2(1+\epsilon )||\sum a_i
t_{k_i}||_T$$
for any sequence of scalars $(a_i)$.\endproclaim

\demo{Proof}  Fix $\alpha <\omega_1$.  Let $(v_i)$ be a normalized
block basic
sequence of $\left(x^{\alpha}_n\right)$ which is
$(1+\epsilon )$-equivalent to the usual unit vector
basis of $c_0$ and let $u_i=v_i/||v_i||$, $i=1,2,\ldots$ .
The idea is to show that
for any sequence of scalars $(a_i)$
$$||\sum a_it_i||_{T,m}\le ||\sum a_iu_i||\tag **$$
and
$$||\sum a_iu_i||_m\le 2(1+\epsilon )||\sum a_it_{2i}||_T\tag ***$$
$m=0,1,\ldots$ , where
$$||\sum a_it_i||_{T,m+1}=\max\left\{||\sum a_it_i||_{T,m},2^{-1}\max_{
\{p_i\}\in \Cal F_1}\sum^k_{i=1}||\sum^{p_{i+1}}_{n=p_i+1}a_nt_n||_{
T,m}\right\},$$
and
$$||\sum a_it_i||_{T,0}=\sup |a_i|.$$
Once this is accomplished we use the fact that there is a constant $
K$,
such that
$$||\sum a_it_i||_T\le ||\sum a_it_{2i}||_T\le K||\sum a_it_i||_T$$
and hence $[u_i]$ is isomorphic to $T$.  (See [C-S,p.26].  In fact the
argument given here is derived from the arguments of Casazza,
Johnson and Tzafriri, [C-S,p.34-38].)

We will establish $(^{**})$ by induction on $m$.  For $m=0$, $(^{*
*})$ is
immediate.  Now assume the inequality holds for $m$ and we will
prove it for $m+1$.

Let $k_n$ be the first element in the support of $u_n$.  Then
$$\align\text{$||\sum a_it_i||_{T,m+1}$}&=\max\left\{||\sum a_it_i
||_{T,m},2^{-1}\max_{\{p_i\}\in \Cal F_1}\sum^k_{i=1}||\sum^{p_{i+
1}-1}_{n=p_i}a_nt_n||_{T,m}\right\},\\
&\le\max\left\{||\sum a_iu_i||,2^{-1}\max_{\{p_i\}\in \Cal F_1}\sum^
k_{i=1}||\sum^{p_{i+1}-1}_{n=p_i}a_nu_n||\right\},\\
&\le\max\left\{||\sum a_iu_i||,2^{-1}\max_{\{p_i\}\in \Cal F_1}\sum^
k_{i=1}||P_i\sum a_nu_n||\right\}\\
&=||\sum a_iu_i||\endalign$$
where $P_i\sum b_nt^{\alpha}_n=$ $\displaystyle{\sum^{q_{i+1}-1}_{
n=q_i}b_nt^{\alpha}_n}$.  The last inequality
above holds because $\{q_i\}=\{k_{p_i}\}\in \Cal F_1$ if $\{p_i\}$ $
\in \Cal F_1$.

For the inequality $(^{***})$ we need to work a little harder.

$$\text{$||\sum a_nu_n||_{m+1}$}=\max\left\{||\sum a_iu_i||_m,2^{-
1}\max_{\{q_i\}\in \Cal F_1}\sum^k_{i=1}||P_i\sum a_nu_n||_m\right
\}$$
where $P_i$ denotes the basis projection onto $[t_j:q_i\le j<q_{i+
1}$].  Fix
$\{q_i\}\in \Cal F_1$ and for each $i$ let $G_i=\left\{n:q_i\le k_
n<k_{n+1},\text{ and }n\notin G_i\right\}$.  For
each $n$ let
$$H_n=\left\{i:k_n\le q_i<k_{n+1}\text{ or }k_n<q_{i+1}\le k_{n+1}
,\text{ and }n\notin G_i\right\}.$$
Consider the sum corresponding to the $q_i$ $^{\prime}s$.
$$ $$
$$\align 2^{-1}\sum^k_{i=1}||P_i\sum a_nu_n||_m&\le 2^{-1}\sum^k_{
i=1}||\sum_{n\in G_i}a_nu_n||_m+2^{-1}\sum^{\infty}_{n=1}\sum_{i\in
H_n}||P_ia_nu_n||_m\\
&\le 2^{-1}\sum^k_{i=1}||\sum_{n\in G_i}a_nu_n||_m+2^{-1}\sum_{n:H_
n\neq\emptyset}2||a_nu_n||_{m+1}\\
&\le 2^{-1}\sum^k_{i=1}2(1+\epsilon )||\sum_{n\in G_i}a_nt_{2k_{n+
1}}||_T+\sum_{n:H_n\neq\emptyset}||a_nt_{2k_{n+1}}||_T\\
&\le (1+\epsilon )\left[\sum^k_{i=1}||\sum_{n\in G_i}a_nt_{2k_{n+1}}
||_T+\sum_{n:H_n\neq\emptyset}||a_nt_{2k_{n+1}}||_T\right]\endalign$$
Observe that there are at most $k$ integers $n$ such that $H_n\neq
\emptyset$ and
that $k\le q_1<k_{m+1}$, where $m$ is the smallest integer such that
$P_1u_m\neq 0$.  Let $n_i=\min G_i$ for $i=1,2,\ldots ,k$.  Then
$$\left\{2k_{n_i+1}:i=1,2,\ldots ,k\right\}\cup\left\{2k_{n+1}:H_n
\neq\emptyset\right\}$$
is a set of at most $2k$ integers greater than $2k$.  Hence this sum in
brackets is at most $2||\sum a_nt_{2k_{n+1}}||_T$.  Therefore
$$\align&\max\left\{||\sum a_iu_i||_m,2^{-1}\max_{\{q_i\}\in \Cal F_
1}\sum^k_{i=1}||P_i\sum a_nu_n||_m\right\}\\
\le&\max\left\{||\sum a_it_{2k_{i+1}}||_T,2(1+\epsilon )\max_{\{q_
i\}\in \Cal F_1}\sum^k_{i=1}||P_i\sum a_nt_{2k_{n+1}}||_T\right\}\\
=&2(1+\epsilon )||\sum a_it_{2k_{i+1}}||_T,\endalign$$
as claimed.\qed\enddemo

\remark{Remark 5.3}  Argyros [A] has modified the construction
to obtain spaces $X_{\alpha}$, $\alpha <\omega_1$, such that
all of the subspaces of $
X_{\alpha}$
have index at least $\alpha$.  To accomplish this he uses the sets $
\Cal F_{\alpha}$ in
the definition of the norm instead of starting with the
space $[x^\alpha_n].$\endremark

\head 6.  Comparison with the averaging index\endhead

In [A-O] the averaging index (See Section 1 for the definition.) was
used
to get a somewhat more constructive version of Mazur's Theorem.
In this section we will show that the averaging index is much larger
than the spreading model index by showing that there exists a Banach
space $X$ such that for every $\alpha <\omega_1$ there is a
weakly null sequence
$(x_n)$ in $X$ with averaging index at least $\alpha$, and yet $\ell^
1$ is not
isomorphic to a subspace of $X$.  We will also show that the spreading
model index can be used to strengthen some of the results in [A-O].

The construction of the example will be based on the infinite
branching James tree construction [J] in combination with $\omega^{
\omega}$.  Let
$$\Cal T=\cup^{\infty}_{n=1}\left\{\left(\alpha_1,\alpha_2,\ldots
,\alpha_n\right):\alpha_i<\omega^{\omega}\text{ for each }i\right\}
.$$
Let $X_0$ be the linear subspace of the functions from $\Cal T$ into $
C_0(\omega^{\omega})$
which are nonzero at only finitely many points of $\Cal T$.  We will
use
the notation $\left(f_{\overline t}\right)$ where $f_{\overline t}
\in C(\omega^{\omega})$ for all $\overline t$ in $\Cal T$ to denote an
element of $X_0$ with the understanding that the index $\overline
t$ runs over $\Cal T$.
We will also use $f_{\overline t}$ to denote the
element of $X_0$ which is 0 except at
$\overline t$ and $f_{\overline t}$ at $\overline t$.

Next we will introduce some linear functionals on $X_0$.  For each
$\overline
s\in\overline {\Cal T}$
and $i<|\overline s|$ (the length of $\overline {s}$) define

$$L_{(\overline s,i)}(f_{\overline t})=\cases 0&\text{if }\overline
s\not
>\overline t\text{ or }|\overline t|<i\\
f_{\overline t}(\overline s(|\overline t|+1))&\text{if }\overline
s>\overline t\text{ and }i\le |\overline t|<|\overline s|\endcases
,$$
where $\overline s>\overline t$ denotes that $\overline s$ is below
$\overline
t$, and extend linearly to $X_0$.

We will refer to a pair $(\overline s,i)$ as a segment and define it to
be the set
of nodes of length at least $i$ which are above $\overline s$.  We also
want to
have a notion of incomparable segments.  Suppose that $j\le i$ and that
$(\overline s,i)$ and $(\overline t,j)$ are segments then $(\overline
t,j)$ and $(\overline s,i)$ are incomparable if
they are disjoint and $\overline t(m)\neq\overline s(m)$ for some $
m\le j$, i.e., they are on
branches that split by level $j$.  Note that if $\left\{\left(\overline
s_k,i_k\right)\right\}$ are pairwise
incomparable nodes then $L_{(\overline s_k,i_k)}(f_{\overline t})\neq
0$ for at most one $k$ for each
$\overline t\in \Cal T$.

Now we will introduce a norm on $X_0$.  For $F\in X_0$ define
$$||F||=\sup\big\{\left[\sum_j\left[L_{(\overline s_j,i_j)}(F)\right
]^2\right]^{\frac 12}:\left\{(\overline s_j,i_j)\right\}\text{ pairwise
incomparable}\big
]$$
Let $X$ be the completion of $X_0$ under this norm.  Clearly $C_0(
\omega^{\omega})$ is
isometric to $X_{\overline t}=\left[f_{\overline t}:f\in C_0(\omega^{
\omega})\right]$ for each fixed $\overline t$.  For each $k\in \Bbb N$
define a projection $P_k$ on $X$ by
$$P_k((f_{\overline
t}))_{\overline s}=\cases f_{\overline s}&\text{if }|\overline s|\le
k\\0&\text{if }|\overline s|>k\endcases.$$
It is easy to see that $||P_k||=1$ for all $k$.

\proclaim{Proposition 6.1}  $\ell^1$ is not isomorphic to a subspace of
$
X$.\endproclaim

The example is similar to an example of Odell [O] and his arguments
can be modified to prove Propositon 6.1.  However we will give a
slightly different proof which does not directly use the branch
functionals.

\demo{Proof}  Suppose that $(y_n)$ is a normalized sequence in $X$
which
is $K$-equivalent to the usual unit vector basis of $\ell^1$.  Because
$
\ell^1$ is
not isomorphic to a subspace of $C_0(\omega^{\omega})$ and range $
(P_k-P_{k-1})$ is
isometric to $\left[\sum C_0(\omega^{\omega})\right]_{\ell^2}$, it can
be shown by induction that $
\ell^1$ is
not isomorphic to a subspace of range $P_k$ for any $k$.  Therefore by
passing to a subsequence we may assume that $(P_k(y_n))_n$ is
weakly Cauchy for each $k$.  Because for each $k$ there are convex
combinations of $(P_k(y_{2n}-y_{2n-1}))$ with small norm we can find a
sequence of disjointly supported (relative to the $y_n$ $^{\prime}
s$)  convex
combinations of $(y_{2n}-y_{2n-1})$, $(z_j)$, and an increasing
sequence of
intergers $(k_j)$ such that
$$\sum^{\infty}_{j=m}||P_{k_m}z_j||<2^{-m},$$
for $m=1,2,\ldots$ .  Moreover because $X_0$ is dense in $X$ we may
assume
(by passing to a subsequence) that
$$\sum^{m-1}_{j=1}||(I-P_{k_m})z_j||<2^{-m},$$
for $m=1,2,\ldots$ .  In this way we get a sequence equivalent to the
unit vector basis of $\ell^1$ which is essentially supported on
disjoint
levels of $\Cal T$.  By a standard perturbation argument we may assume
that $(I-P_{k_{m+1}})z_m=0=P_{k_m}z_m$, and that $z_{m|_{\overline
t}}\neq 0$ for only finitely
many $\overline t$ for all $m$.

Next note that by a theorem of James [J] we may assume that $(z_j)$
is $(1+\epsilon )$-equivalent to the basis of $\ell^1$.  For any node
$\overline
s$ let
$\Cal W(\overline s)=\left\{\overline t:\overline t>\overline s\right
\}$, the wedge determined by $\overline s$.  For each $i$ there are
finitely many nodes $\overline s(i,j)$ of length $k_i$ such that if
$\displaystyle{\overline t\notin\cup_j\Cal W(\overline s(i,j))}$ then $
z_{i|_{\overline t}}=0$.  Let $N$ be the number of
nodes in the support of $z_1$.  We claim that for each $i$ there is a
set
$\Cal F=\Cal F(i)$ of cardinality at most $N$ such that
$$||z_i-z_i|_{\cup_{j\in \Cal F}\Cal W(\overline s(i,j))}||<4\epsilon
.$$
Indeed, if not let $\Cal S=\left\{\left(\overline s,i\right)\right
\}$ be a family of incomparable segments
which compute the norm of $z_1+z_i$.  Let $\Cal S^{\prime}$ be the set
of segments in
$\Cal S$ which intersect the support of $z_1$ and $\Cal
S^{\prime\prime}
=\Cal S\backslash \Cal S^{\prime}$.  Clearly $\Cal S^{\prime}$ contains
at most $N$ segments.  We have that
$$4(1+\epsilon )^{-2}\le ||z_1+z_i||^2$$
$$=\sum_{(\overline s,j)\in \Cal S^{\prime}}\left[L_{(\overline
s,j)}(z_1+z_i)\right]^2+\sum_{(\overline s,j)\in \Cal S^{
\prime\prime}}\left[L_{(\overline s,j)}(z_1+z_i)\right]^2$$
$$\le\left[\left[\sum_{(\overline s,j)\in \Cal S^{\prime}}\left
[L_{(\overline s,j)}(z_1)\right]^2\right]^{\frac
12}+\left[\sum_{(\overline
s,j)\in \Cal S^{\prime}}\left[L_{(\overline s,j)}(z_i)\right
]^2\right]^{\frac 12}\right]^2$$
$$+\sum_{(\overline s,j)\in \Cal S^{\prime\prime}}\left[
L_{(\overline s,j)}(z_1+z_i)\right]^2$$
$$\le (1+(1-4\epsilon ))^2+(4\epsilon )^2=4(1-4\epsilon +8\epsilon^
2).$$
Clearly this is impossible for small enough $\epsilon$.

It follows that by another perturbation argument that we may
assume that each $z_i$ is supported in at most $N$ wedges.  As above
let
$\overline s(i,j),\text{ }j=1,2,\ldots ,N$ be the nodes of length $k_i$
so that $
z_i$ is supported
in $\displaystyle{\cup_j\Cal W(\overline s(i,j))}$.  We will next
refine our sequence $
(z_i)$
to get a subsequence such that there are branches $\overline b_1,\ldots
,\overline b_k$, $k\le N$
such that if $\overline s$ is any branch then $\overline s>\overline
s(i,j)$ for at most $N$ nodes not
on some $\overline b_m$ and $\overline b_j>\overline s(i,j)$ for $
j=1,2,\ldots ,k$ for all $i$.  Such a
subsequence is easily determined by induction on $N$.  Indeed, if $
N=1$,
either there are infinitely many $i$ and incomparable branches
$(\overline
t_i,k_i)$
such that $\overline t_i>\overline s(i,1)$ and $\overline t_i\not
>\overline s(m,1)$ for any $m\neq i$, or there is a
branch $\overline b_1$ which contains all but finitely many of the
nodes $\overline
s(i,1)$.
Now suppose that $\overline b_1,\ldots ,\overline b_k$ are branches
such that if $\overline
s$ is any branch
then $\overline s>\overline s(i,j)$ for at most N-1 nodes not on some
$\overline
b_m,$ and $\overline b_j>\overline s(i,j)$
for $j=1,2,\ldots ,k$ for all $i$.  As above if there is some branch
which
contains all but finitely many of the nodes $\overline s(i,N)$ we add
that branch
to our list as $\overline b_{k+1}$ and pass to a subsequence $(z_i
)_{i\in M}$ such that
$\overline b_{k+1}>\overline s(i,N)$ for all $i\in M$.  If this is not
the case then we can find
a subsequence such that there is at most one of the nodes $\overline
s(i,N)$ on
any branch.  Clearly this subsequence has the required properties.

To complete the argument we need to make a few observations about
the norm on $X$.  First observe that if $\overline s$ is any branch and
for each
node $\overline t$ on $\overline s$, $f_{\overline t}$  is a fixed
function in \linebreak $
C_0(\omega$ $^{\omega})$ then $[f_{\overline t}]_{\overline s>\overline
t}$ is
isomorphic to $c_0$.  Second if the norm of $\sum a_iz_i$ is computed
using
only segments which intersect at most $N$ of the supports of the $
z_i$ $^{\prime}s$,
i.e., $\overline s(i,j)$ for at most $N\text{ }j^{\prime}s$, then the
norm is at most $\left
(\sum a^2_i\right)^{\frac 12}N^{\frac 12}$.
Finally observe that for the sequence $(z_i)$ a segment can intersect
more than $N$ of the supports of the $z_i$ $^{\prime}s$ only if it lies
along one of
the branches $\overline b_k$.  A straightforward computation using
these
observations shows that the $z_i$ $^{\prime}s$ are not equivalent to
the unit
vector basis of $\ell^1$.\qed\enddemo

Our next goal is to show that the averaging index of $X$ is
uncountable.  The basic idea is to construct for each $\alpha <\omega_
1$ a weakly
null sequence with averaging index $\alpha$ by using well-founded
subtrees
of $\Cal T$ of order $\alpha$.

\proclaim{Proposition 6.2}  Let $\Cal S$ be a well-founded subtree of $
\Cal T$ and
for each $\overline s\in \Cal S$ let $\left[f^{^{n}}_{\overline
s}\right]$ be a normalized weakly null sequence in $
X\frac {}s$.
Then $\left[f^{^{n}}_{\overline s}\right]_{n\in \Bbb N,\overline s\in
\Cal S}$ (reordered) is a weakly null sequence in $
X$.\endproclaim

\demo{Proof}  $X^{*}$ is the closed linear span of the functionals
$\left[L_{(\overline t,i)}\right]_{\overline t\in\overline {\Cal T}
,i\in \Bbb N}$.  Because $\Cal S$ is well founded, for any $\overline
t\in\overline {\Cal T}$, $\overline s<\overline t$ for
only finitely many $\overline s\in \Cal S$.  Hence for any $\epsilon
>0$, $\overline t\in\overline {\Cal T}$, and $i\in \Bbb N$,
$|L_{(\overline t,i)}\left[f^{^{n}}_{\overline s}\right]|\ge\epsilon$
for only finitely many $
n$ and $\overline s$.  Therefore
$\left[f^{^{n}}_{\overline s}\right]_{n\in \Bbb N,\overline s\in \Cal
S}$ is a weakly null sequence in $
X$.\qed

In order to estimate the averaging index we need to have some
information about the $w^{*}$ topology on the functionals
$L_{(\overline
t,i)}$.\enddemo

\proclaim{Lemma 6.3}  Let $\overline t\in \Cal T$ and let $(\alpha_
i)$ be a sequence in $\omega^{\omega}$
which converges to $\alpha$.  Then
\itemitem{i)}  If $\alpha <\omega^{\omega},$ $L_{(\overline t+(\alpha_
i),j)}\overset {w^{*}}\to {\longrightarrow}$ $L_{(\overline t+(\alpha
),j)}$
\itemitem{ii)} If $\alpha =\omega^{\omega},$ $L_{(\overline t+(\alpha_
i),j)}\overset {w^{*}}\to {\longrightarrow}L_{(\overline t,j)}$.
\endproclaim

\demo{Proof}  We need only consider the values of the functionals at
$f_{\overline s}$ for those $\overline s\in \Cal T$ such that
$\overline
s\le\overline t$ and $|\overline t|-1\le |\overline s|\le |\overline
t|$.  For all others
the values do not depend on $(\alpha_i)$.  If $\overline s=\overline
t$ then
$$L_{(\overline t+(\alpha_i),j)}\left[f_{\overline
s}\right]=f_{\overline s}
(\alpha_i)\longrightarrow f_{\overline s}(\alpha )=L_{(\overline
t+(\alpha
),j)}\left[f_{\overline s}\right],\text{ if
}\alpha_i\longrightarrow\alpha
<\omega^{\omega},$$
and
$$L_{(\overline t+(\alpha_i),j)}\left[f_{\overline
s}\right]=f_{\overline s}
(\alpha_i)\longrightarrow f_{\overline s}(\alpha )=0=L_{(\overline
t,j)}\left
[f_{\overline s}\right],\text{ if
}\alpha_i\longrightarrow\omega^{\omega}
.$$
If $\overline t=\overline s+(\beta )$, then
$$L_{(\overline t+(\alpha_i),j)}\left[f_{\overline
s}\right]=f_{\overline s}
(\beta )\longrightarrow f_{\overline s}(\beta )=L_{(\overline t+(\alpha
),j)}\left[f_{\overline s}\right],\text{ if
}\alpha_i\longrightarrow\alpha
<\omega^{\omega},$$
and
$$L_{(\overline t+(\alpha_i),j)}\left[f_{\overline
s}\right]=f_{\overline s}
(\beta )\longrightarrow f_{\overline s}(\beta )=L_{(\overline
t,j)}\left
[f_{\overline s}\right],\text{ if
}\alpha_i\longrightarrow\omega^{\omega}
.$$\qed\enddemo

We will need to use well-founded subtrees of $\Cal T$ of a special
type, so
as a technical convenience we introduce the following.

\definition{Definition}  A well-founded tree $\Cal S\subset \Cal T$ is
said to be complete
if
\itemitem{i)}  $\overline s+(\beta )\in \Cal S$ for some $\beta
<\omega^{
\omega}$ and $\overline s\in \Cal T$ then $\overline s+(\alpha )\in
\Cal S$ for
all $\alpha <\omega^{\omega}$ and $\overline s\in \Cal S$,
\itemitem{ii)}  $(\alpha )\in \Cal S$ for all $\alpha
<\omega^{\omega}$.
\vskip .1in
Next we will verify that such things exist.\enddefinition

\proclaim{Proposition 6.4}  For every $\alpha <\omega_1$ there exists a
complete
well-founded tree $\Cal S\subset \Cal T$ of order at least
$\alpha$.\endproclaim

\demo{Proof}  For $\alpha =1$ this is obvious.  Suppose for all $
\beta <\alpha$ there
is a complete well-founded subtree $\Cal S$ of $\Cal T$ with order $
\beta$.  If $\alpha =\beta +1$,
we define $\Cal U=\left\{(\eta )+\overline s:\overline s\in \Cal S
\cup \{()\} \right.$ and \linebreak $ \left.\eta <\omega^{\omega}\right\}$. 
 Clearly $o
(\Cal U)=o(\Cal S)+1$
and it is obvious that $\Cal U$ is complete.  If $\alpha$ is a limit
ordinal and
$\alpha_i\uparrow\alpha$, for each $i$ let $\Cal S_{\alpha_i}$ be a
complete well-founded subtree of $
\Cal T$ of
order $\alpha_i$.  Let $\Cal U=\left\{(\eta )+\overline s:\overline
s\in\cup_i\Cal S_{\alpha_i}\cup \{()\}\right.\text{ and }
\left.\eta <\omega^{\omega}\right\}$.
Clearly $\Cal U$ is
complete, $o(\Cal U$ $)\ge o(\Cal S_{\alpha_i})$ for all $i$ and thus $
o(\Cal U)\ge\alpha$.\qed\enddemo

We are now ready to show that there are weakly null sequences in
$X$ with large averaging index.

\proclaim{Proposition 6.5}  Suppose that $\Cal S$ is a complete
well-founded
subtree of $\Cal T$ and for each $\overline s\in \Cal S,$ $\left[f^
{^{n}}_{\overline s}\right]_{n\in \Bbb N}$ is the Schreier sequence.
Then for every $\beta <o(\Cal S)$ and $\overline s\in \Cal S^{(\beta
)}$, $L_{(\overline s,1)}\in$\linebreak $A_{\beta}\left[\frac
12,\left[f^
{^{n}}_{\overline t}\right]_{n\in \Bbb N,\overline t\in \Cal
S},B_{X^{*}}\right]$.\endproclaim

\demo{Proof}  We proceed by induction on $\beta$.  It is sufficient to
prove the result for $\beta +1$ assuming it for $\beta$.  If $\overline
s\in \Cal S^{(\beta +1)}$, then
because $\Cal S$ is complete $\overline s+(\alpha )\in \Cal S^{(\beta
)}$ for all $\alpha <\omega^{\omega}$.  Moreover
$L_{(\overline s+(\alpha ),1)}\left[f^{^{n}}_{\overline
s}\right]=f^{^{n}}_{\overline s}(\alpha
)$.  In Section 4 it was shown that
$\omega^{\omega}\in \Cal O^1\left(\frac 12,\left(f^n\right),\omega^{
\omega}\right)\subset A_1\left(\frac
12,\left(f^{^{n}}\right),\omega^{\omega}\right
)$ where $\left(f^{^{n}}\right)$ is the Schreier
sequence.  Therefore by Lemma 6.3,
$L_{(\overline s,1)}\in A_1\left[\frac 12,\left[f^{^{n}}_{\overline
s}\right]_{
n\in \Bbb N},A_{\beta}\left[\frac 12,\left[f^{^{n}}_{\overline
t}\right]_{n\in
\Bbb N,\overline t\in \Cal S},B_{X^{*}}\right]\right]$.  Because the
definition of the averaging index only requires that a subsequence
have $\ell^1$ $-SP\ge\frac 12$, it follows that
$$\align L_{(\overline s,1)}&\in A_1\big[\frac
12,\left[f^{^{n}}_{\overline t}\right
]_{n\in \Bbb N,\overline t\in \Cal S},A_{\beta}\left[\frac 12,\left
[f^{^{n}}_{\overline t}\right]_{n\in \Bbb N,\overline t\in \Cal
S},B_{X^{*}}\right]\\
&=A_{\beta +1}\left[\frac 12,\left[f^{^{n}}_{\overline t}\right]_{n\in
\Bbb N,\overline
t\in \Cal S},B_{X^{*}}\right].\qed\endalign$$\enddemo

\proclaim{Corollary 6.6}  For every $\alpha <\omega_1$ there is a
weakly null
sequence in $X$ such that \linebreak $o(A,\frac 12$
$)\le\alpha$.\endproclaim
\vskip .1in
Next we want to consider the relationship between the spreading
model index and the constructive version of Mazur's theorem as
presented in [A-O].  One purpose of that paper was to try to use the
averaging index to determine if given a Banach space $X$ there is an
integer $k$ such that any weakly null sequence needs to be averaged at
most $k$ times in order to get a norm null sequence.  In what follows
we will use the notation of [A-O] and refer the reader there for the
relevant definitions.

The definition of the spreading model index makes some arguments of
the type used in [A-O] difficult because of its sensitivity to passing
to subsequences.  On the other hand we are going to consider all
weakly null sequences in the space $X$ so it seems natural to look for
sequences which are in some sense extremal.

\proclaim{Proposition 6.7}  Suppose $X$ is a subspace of $C(K)$ for
some
compact metric space $K$ and that $(x_n)$ is a weakly null sequence in
the unit ball of $X$.  Then for every $\epsilon >0$ there is a
subsequence
$(x_n)_{n\in M}$ of $(x_n)$ such that for any $L\subset M$, infinite,
and $
\alpha <\omega_1$,
$$S^{\alpha}\left(\epsilon ,\left(x_n\right)_{n\in L},K\right)=S^{
\alpha}\left(\epsilon ,\left(x_n\right)_{n\in M},K\right).$$
\endproclaim

\demo{Proof}  The idea is to use repeatedly the following lemma
[A-O].\enddemo

\proclaim{Lemma 6.8}  Let $K$ be a second countable compact
Hausdorff space and let $(x_n)$ be a weakly null sequence in $C(K)$.
Then
there is a subsequence $(x_n)_{n\in M}$ such that for every $t\in
K$ and every
neighborhood $\Cal N$ of $t$ there is a neighborhood $\Cal N^{\prime}$
of $
t$, $\Cal N^{\prime}$ $\subset \Cal N$, such that
$(x_{n|\Cal N^{\prime}})$ has a spreading model.\endproclaim

We will prove the proposition by induction on $\alpha$.  Note that
there is
an ordinal $\alpha_0$ such that $S^{\alpha_0}\left(\epsilon ,\left
(y_n\right),K\right)=\emptyset$ for every $\epsilon >0$ and
weakly null sequence $(y_n)$ in $C(K)$, and thus the induction is over
a
countable set.  Fix $\epsilon >0$.

\definition{Induction Hypothesis}  If $(x_n)$ is a weakly null sequence
in the unit
ball of $X$.  Then for every $\epsilon >0$ there is a subsequence $
(x_n)_n$ $_{\in M}$ of
$(x_n)$ such that for any $L\subset M$, infinite, and $\beta\le\alpha$,
$$S^{\beta}\left(\epsilon ,\left(x_n\right)_{n\in
L},K\right)=S^{\beta}\left
(\epsilon ,\left(x_n\right)_{n\in M},K\right).$$\enddefinition

Observe that the lemma proves the result for $\alpha =1$.  Now suppose
that the result is true for all $\beta <\alpha$.  If $\alpha =\beta
+1$, for some $\beta$, then
let $(x_n)_{n\in L}$ be a subsequence such that
$$S^{\lambda}\left(\epsilon ,\left(x_n\right)_{n\in J},K\right)=S^{
\lambda}\left(\epsilon ,\left(x_n\right)_{n\in L},K\right),$$
for all infinite $J\subset L$ and $\lambda\le\beta$.  Apply the lemma
to
$\left(x_{n|_{S^{\beta}\left(\epsilon ,\left(x_n\right)_{n\in
J},K\right
)}}\right)_{n\in L}$ to get a subsequence $(x_n)_{n\in M}$ and note
that
the case $\alpha =1$ applies to show that
$$S^{\beta +1}\left(\epsilon ,\left(x_n\right)_{n\in J},K\right)=S^
1\left(\epsilon ,\left(x_{n|_{S^{\beta}\left(\epsilon ,\left(x_n\right
)_{n\in L},K\right)}}\right)_{n\in J},S^{\beta}\left(\epsilon ,\left
(x_n\right)_{n\in L},K\right)\right)$$
$$=S^1\left(\epsilon ,\left(x_{n|_{S^\beta{\left(\epsilon ,\left(x_
n\right)_{n\in L},K\right)}}}\right)_{n\in M},S^{\beta}\left(\epsilon
,\left(x_n\right)_{n\in L},K\right)\right)=S^{\beta +1}\left(\epsilon
,\left(x_n\right)_{n\in M},K\right).$$
for all infinite $J\subset M$.

If $\alpha$ is a limit ordinal, let $\alpha_i\uparrow\alpha$ and choose
infinite sets
$M_1\subset M_2\subset\ldots\subset M_i\subset\ldots$, such that
$$S^{\lambda}\left(\epsilon ,\left(x_n\right)_{n\in J},K\right)=S^{
\lambda}\left(\epsilon ,\left(x_n\right)_{n\in M_i},K\right),$$
for all infinite $J\subset M_i$ and $\lambda\le\alpha_i$.  Then if $
M$ is an infinite set such
that $M\backslash M_i$ is finite for all $i$, it follows that
$$S^{\lambda}\left(\epsilon ,\left(x_n\right)_{n\in J},K\right)=S^{
\lambda}\left(\epsilon ,\left(x_n\right)_{n\in M_i},K\right),$$
for all infinite $J\subset M$ and $\lambda <\alpha$.  However because
$\displaystyle{S^{\alpha}=\cap_{\lambda <\alpha}S^{\lambda}}$, the
equality holds for $ \lambda =\alpha$ as well.\qed

\proclaim{Corollary 6.9}  Suppose that $(x_n)$ is a weakly null
sequence
in the ball of $C(K)$ for some compact metric space $K$.  Then there is
a subsequence $\left(x_n\right)_{n\in M}$ such that
$$A^{\lambda}\left(\epsilon ,\left(x_n\right)_{n\in M},K\right)=S^{
\lambda}\left(\epsilon ,\left(x_n\right)_{n\in M},K\right),$$
for all $\lambda\le\alpha$ and $\epsilon >0$.\endproclaim

\demo{Proof}  The proof above gives us a subsequence so that
passing to subsequences has no effect on $\ell^1-SP\left(x_{n|\Cal
N}\right
)$ for some $\Cal N\downarrow \{t\}$
for all $t$, therefore
$$A^1\left(\epsilon ,\left(x_n\right)_{n\in
M},K\right)=S^1\left(\epsilon
,\left(x_n\right)_{n\in M},K\right).$$
Induction completes the proof.\qed\enddemo

\proclaim{Corollary 6.10}  Let $X$ be a subspace of $C(K)$, $K$ compact
metric, such that\linebreak
$\sup\left\{o\left(S,\left(x_n\right),\epsilon\right):\epsilon >0\text{
and }\left
(x_n\right)\subset B_X\text{ is weakly null}\right\}<\omega^{k+1}$, for
some\linebreak $k\in\left\{-1,0,1,2,\ldots\right\}$.  Then $X$ has
property-$
A(k+2)$.\endproclaim

\demo{Proof}  To establish property $A(m)$ it is sufficient to show
that every weakly null sequence has a subsequence with a norm null
convex block subsequence of $m$-Averages.  According to the previous
corollary every weakly null sequence has a subsequence for which
the averaging index and spreading model index agree.  Therefore by
[A-O, Theorem 4.1], $X$ has property-$A(k+2)$.\qed\enddemo

Because the spreading model index is in general smaller than the
averaging index this corollary gives a real strengthening of [A-O,
Corollary 4.2].  In particular the space constructed above has
averaging index $\omega_1$ but spreading model index at most
$\omega$.\newline
\vskip .13in
Acknowledgement.  The authors were visitors at the University of
Texas during the time that a portion of this work was done and
would like to thank the Department of Mathematics at University of
Texas for its hospitality.  In particular we are grateful to H.P.
Rosenthal and E. Odell for many helpful comments.

\Refs

\item{[A]}  S. Argyros, {\it Banach spaces of the type of Tsirelson\/},
preprint.

\item{[A-O]}  D. Alspach and E. Odell, {\it Averaging weakly null
sequences\/}, Lecture Notes in mathematics Vol. 1332, Springer-Verlag,
Berlin 1988.

\item{[Bo]}  J. Bourgain, {\it On convergent sequences of continuous
functions\/}, Bull. Soc. Math. Bel., 32 (1980), 235-249.

\item{[C-S]}  P. Casazza and T. Shura, ``Tsirelson's Space'', Lecture
Notes in Mathematics Vol. 1363, Springer-Verlag, Berlin/New York 1989.

\item{[D1]}  J. Diestel, ``Geometry of Banach Spaces--Selected
Topics'',
Lecture Notes in Mathematics Vol. 485, Springer-Verlag, Berlin/New
York 1975.

\item{[D2]}  J. Diestel, ``Sequences and Series in Banach Spaces'',
Graduate Texts in Mathematics Vol. 92, Springer-Verlag, Berlin/New
York 1984.

\item{[Do]}  L. Dor, {\it On projections in $L_1$\/}, Ann. of Math.,
Vol. 102
(1975), 463-474.

\item{[G-H]}  D. C. Gillespie and W. A. Hurwitz, {\it On sequences of
continuous functions having continuous limits\/}, Trans. AMS, 32
(1930),
527-543.

\item{[H-O-R]}  R. Haydon, E. Odell, and H. P. Rosenthal, {\it On
certain
classes of Baire-1 functions with applications to Banach space
theory\/}, Lecture Notes in Mathematics Vol. 1470, Springer-Verlag,
Berlin 1991.

\item{[J]}  R. C. James, {\it A separable somewhat reflexive Banach
space with nonseparable dual\/}, Bull. AMS 80 (1974), 738-743.

\item{[K-L]}  A. S. Kechris and A. Louveau, {\it A classification of
Baire-1 Functions\/}, Tran. AMS 318 (1990), 209-236.

\item{[K]}  K. Kuratowski, ``Topology, I'', Academic Press Inc., New
York 1966.

\item{[L-S]}  J. Lindenstrauss and C. Stegall, {\it Examples of
separable
reflexive spaces which do not contain $\ell_1$ and whose duals are
nonseparable,\/} Studia Math. 54 (1975), 81-105.

\item{[L-T,I]}  J. Lindenstrauss and L. Tzafriri, ``Classical Banach
Spaces I:  Sequence Spaces'', Ergibnisse der Mathematik 92,
Springer-Verlag, Berlin/New York 1977.

\item{[L-T,II]}  J. Lindenstrauss and L. Tzafriri, ``Classical Banach
Spaces II:  Function Spaces'', Ergibnisse der Mathematik 97,
Springer-Verlag, Berlin/New York 1979.

\item{[M-R]}  B. Maurey and H. P. Rosenthal, {\it Normalized weakly
null
sequences with no unconditional subsequence,\/} Studia Math. 61 (1977)
77-98.

\item{[O]}  E. Odell, {\it A normalized weakly null sequence with no
shrinking subsequence in a Banach space not containing $\ell_1$, \/}
Compositio Math. 41 (1980), 287-295.

\item{[P-S]}  A. Pelczynski and W. Szlenk, {\it An example of a
non-shrinking basis~\/}, Rev. Roumaine math. pures et appl. 10 (1965),
961-966.

\item{[R]}  H. P. Rosenthal, {\it A characterization of Banach spaces
containing $\ell^1$\/}, Proc. Nat. Acad. Sci. U. S. A. 71 (1974),
2411-2413.

\item{[Sch]}  J. Schreier, {\it Ein Gegenbeispiel zur Theorie der
schwachen Konvergenz\/}, Studia Math. 2 (1930), 58-62.

\item{[Sz]}  W. Szlenk, {\it The non-existence of a separable reflexive
Banach space universal for all separable reflexive Banach spaces, \/}
Studia Math. 30 (1968), 53-61.

\item{[Z]}  Z. Zalcwasser, {\it Sur une propriet\'e der champ des
fonctions continues\/}, Studia Math. 2 (1930), 63-67.

\endRefs

{\smc Department of Mathematics,  Oklahoma State University,
Stillwater, OK 74078-0613

Department of Mathematics, University of Crete,
Heraklion, Crete}
\bye